\documentclass[12pt,reqno]{amsart}
\usepackage[cp1251]{inputenc}
\usepackage[russian]{babel}
\usepackage{amsmath}
\usepackage{amssymb}
\usepackage{amsfonts}
\usepackage{graphicx}
\usepackage{cmap,cite}

\newtheorem{lemma}{Лемма}
\newtheorem{theorem}{Теорема}

\def\klabel#1{}

\begin{document}

\leftline{517.984.50}
\begin{flushright}{Посвящается Сергею Львовичу Соболеву}
\end{flushright}
\begin{center}
 {\Large Пространства Соболева
и операторы вихрь и градиент дивергенции }
  \end{center}
\begin{center}
\end{center}
\begin{center}
Р.~С.~Сакс
\end{center}
\begin{center}
Институт математики c ВЦ УФИЦ РАН,\\
ул.Чернышевского, 112,
450077, г. Уфа, Россия
\end{center}
\begin{center}
romen-saks@yandex.ru
\end{center}
{\bf Аннотация. }{\small 
	 В ограниченной области $G\subset \textrm{R}^3$ с гладкой
границей изучаются  краевые и спектральные задачи для операторов ротор (вихрь) и градиент дивергенции   $+\lambda\,I$  в пространствах Соболева. 

 При $\lambda\neq 0$ операторы  расширяются 
 (методом Б.Вайнберга и В.Грушина)   до эллиптических матриц,
а краевые задачи удовлетворяют условиям эллиптичености В.Солонникова.
Из теории и оценок вытекают полезнные свойства решений   спектральных задач.
	Операторы $\nabla \text{div}$ и $ \text{rot}$  имеют самосопряженные
расширения $\mathcal{N}_d$ и $\mathcal{S}$ в  ортогональные
подпространства $\mathcal{A}_{\gamma }$ и $\mathbf{V}^0$  потенциальных и вихревых  полей в  $\mathbf{L}_{2}(G)$, а их собственные векторы задают ортогональные базисы в 
 $\mathcal{A}_{\gamma }$ и $\mathbf{V}^0$,  элементы которых представляются   
    рядами Фурье,  а операторы -  преобразованиями рядов.
Определены  аналоги пространств Соболева $\mathbf{A}^{2k}_{\gamma }$ и $\mathbf{W}^m$ 
  порядков  $2k$ и  $m$  в классах  потенциальных и вихревых  полей
 и  классы  $ C(2k,m)$ их прямых сумм.
 Доказано, что   при  $\lambda\neq Sp(\mathrm{rot})$ оператор  $ \text{rot}+\lambda\,I$  отображает  класс  $C(2k,m+1)$ на класс $C(2k,m)$
 взаимно однозначно и непрерывно; а   при  $\lambda\neq Sp(\nabla \mathrm{div})$ оператор $\nabla \text{div}+\lambda\,I$ отображает  класс  $C(2(k+1), m)$ на класс $C(2k,m)$,  соответственно.}
 
 {\small  Ключевые слова:  пространства  Соболева, операторы градиент, дивергенция,  ротор,
 	эллиптические  краевые и спектральные задачи.}

 {\bf Abstract }  Saks R.S.
  
 Sobolev spaces and operators vorticity and the gradient of the divergence
 
{\small   In a bounded domain $G$ with smooth
 border studied 
 boundary value and spectral problems for operators of the rotor (vortex) and the gradient of the divergence  $+\lambda\,I$ in the Sobolev spaces. 
 For $\lambda\neq 0$ these operators  are
reducible  ( by B. Veinberg and V. Grushin method)
 to elliptical matrices
 and the boundary value problems satisfy the conditions of V. Solonnikov's ellipticity.
 Useful properties of solutions of these spectral problems follow from the theory and estimates.
     The $\nabla \text{div}$ and $ \text{rot}$ operators have self-adjoint
 extensions $\mathcal{N}_d$ and $\mathcal{S}$ in orthogonal
 subspaces $\mathcal{A}_{\gamma }$ and $\mathbf{V}^0$  which formed from potential and vortex fields in $\mathbf{L}_{2}(G)$.
 Their eigenvectors forme orthogonal
 basis in  $\mathcal{A}_{\gamma }$ and $\mathbf{V}^0$
 elements of which are presented   by
 Fourier series  and operators are transformations of series. 
 We define
analogues of Sobolev spaces $\mathbf{A}^{2k}_{\gamma }$ and $\mathbf{W}^m$ 
orders of $2k$ and $m$   in classes of potential and vortex fields 
and classes $ C (2k,m)$ of their direct sums.
It is proved that if $\lambda\neq Sp(\mathrm{rot})$ the operator $ \text{rot}+\lambda\,I$ displays the class $C(2k,m+1)$ on the class $C(2k,m)$ one-to-one and continuously.  And  if  $\lambda\neq Sp(\nabla \mathrm{div})$ operator $\nabla \text{div}+\lambda\,I$  maps class  $C(2(k+1), m)$  on the class $C(2k,m)$,  respectivly.}

{\small Keywords: the Sobolev spaces, 
	operators  gradient,  divergence and curl (rotor), elliptic
	boundary value problems, spectral problems.}

  {\small Эта работа
	есть продолжение исследований автора\,\cite{saUMJ13}, \cite{sa15a}-\cite{ds18}.(см. также arXiv:1704.05699v1,\,
	arXiv:1710.06428v1 \,и \,arXiv:1712.03804).}
\section { Основные подпросранства $\mathbf{L}_{2}(G)$} 
\subsection {Шкала  пространств Соболева} Мы рассматриваем линейные
пространства над полем $\mathbb{C}$ комплексных чисел. Через
$\mathbf{L}_{2}(G)$ обозначаем пространство Лебега вектор-функций  (полей),
квадратично интегрируемых в $G$ с внутренним произведением
 $(\mathbf {u},\mathbf {v})= \int_G \mathbf
{u}\cdot\overline{\mathbf {v}} d \mathbf {x}$ и нормой
$\|\mathbf{u}\|~= (\mathbf {u},\mathbf {u})^{1/2}$.
 Пространство Соболева, состоящее из полей,
 принадлежащих $\mathbf{L}_{2}(G)$ вместе с обобщенными производными
 до порядка $ s> 0$, обозначается через
$\mathbf{H}^{s}(G)$, $\|\mathbf {f}\|_s$ -норма его элемента
$\mathbf {f}$; $\mathbf{H}^{0}(G)\equiv\mathbf{L}_{2}(G)$. Замыкание
в $\mathbf{H}^{s}(G)$ множества $\mathcal{C}^{\infty}_0(G)$
обозначается через $\mathbf{H}^{s}_0(G)$. Пространство Соболева
отрицательного порядка $\mathbf{H}^{-s}(G)$ двойственно к
$\mathbf{H}^{s}_0(G)$
 (см. пространство $W_p^{(l)}(\Omega)$ при $p=2$ в $\S 3$ гл. 4 \cite{sob},
 $H^k(Q)$ в $\S 4$ гл. 3 \cite{mi}).

 В области $G$ с гладкой границей $\Gamma$ в каждой точке $y\in\Gamma$
 определена нормаль $\mathbf {n}(y)$ к $\Gamma$.
Поле $\mathbf {u}$ из $\mathbf{H}^{s+1}(G)$ имеет след $
\gamma(\mathbf {n}\cdot\mathbf {u})$ на $\Gamma$ его нормальной
компоненты, который принадлежит пространству Соболева-Слободецкого
$\mathbf{H}^{s+1/2}(G)$, $|\gamma(\mathbf {n}\cdot\mathbf
{u})|_{s+1/2}$- его норма.

\subsection {Разложение $\mathbf{L}_{2}(G)$ на два класса  $\mathcal{A}$ и  $\mathcal{B}$  потенциальных и соленоидальных полей. }
Пусть $h$- функция  из ${H}^{1}(G)$, а $\mathbf{u}=\nabla h$ - ее
градиент.  Обозначим ${\mathcal{{A}}}(G) =\{\nabla h, h\in H^1(G)\}$
-- подпространство в $\mathbf{L}_{2}(G)$, а через ${\mathcal{{B}}}(G)$ - его ортогональное дополнение.  Соотношения
 $(\mathbf {u},\nabla h)=0$ для любой   $ h\in H^1(G)$ означают, что $\mathcal{B}(G)=\{\mathbf{u}\in\mathbf{L}_{2} (G): \mathrm{div} \mathbf{u}=0, \,\gamma(\mathbf{n}\cdot \mathbf{u})=0 \}.$ Итак, 
\begin{equation*}\klabel{wor 1}\mathbf{L}_{2}(G)=
{\mathcal{{A}}}(G)\oplus{\mathcal{{B}}}(G).
 \footnote {{\small Это разложение я взял из статьи   Z.Yoshida и Y.Giga  \cite{yogi}.  Авторы5 называют его разложением Вейля \cite{hw}, а  ${\mathcal{{B}}}(\Omega)$  обозначают как ${L}_{\sigma}^2(\Omega)$.}}
  \eqno{(1.1)}\end{equation*}
 {\small В разложении Г.Вейля:   ${L}_2(G)\equiv\mathfrak{F}_0=
 \mathfrak{G}+\mathfrak{F}'$, где
  $\mathfrak{G}$ есть
 замыкание в норме ${L}_2$ градиентов $\nabla\psi$ функций $\psi\in
\mathcal{C}_0^1(G)$, а $\mathfrak{F}'$-множество соленоидальных
элементов в $\mathfrak{F}_0$ (\cite{hw} Теорема II)}.

 Если граница
области $G$ имеет положительный род $\rho$, то $ \mathcal{B}$
содержит в себе конечномерное подпространство
\begin{equation*}\klabel{bh 1} \mathcal{B}_H=\{\mathbf{u}\in\mathbf{L}_{2}
(G): \mathrm{rot} \mathbf{u}=0, \,\mathrm{div} \mathbf{u}=0,
 \,\gamma(\mathbf{n}\cdot \mathbf{u})=0 \}. \end{equation*} Его
размерность
 равна $\rho$\cite{boso}, а базисные поля
 $\mathbf{h}_j\in \mathcal{C}^\infty(G)$ \cite{hw}.

 Ортогональное дополнение в $\mathcal{B}$  к   $\mathcal{B}_H$   назовем классом вихревых полей и  
 обозначим  $\mathbf {V}^{0} (G)$. Значит,
\begin{equation*}\klabel{bhr 1}\mathcal{B}(G)=\mathcal{B}_{H} (G)\oplus
\mathbf {V} ^{0} (G). \footnote {В\cite{yogi}:
${L}_{\sigma}^2(\Omega)={L}_{\Sigma}^2(\Omega)\oplus{L}_{H}^2(\Omega)$. Символ $L$ перегружен. Мы изменили авторские обозначения пространств
${L}_{\Sigma}^2(\Omega)$
 и ${L}_{H}^2(\Omega)$ на
 $\mathbf {V}^{0} (\Omega)$ и $\mathcal{B}_{H} (\Omega)$. } \eqno(1.2) \end{equation*}

 По определению 
 $\mathcal{A}_{\gamma} = \{\nabla h, h\in H^2(G), \gamma(\mathbf {n}\cdot {\nabla}) h=0 \}.$
 
  {\small  С.Л.Соболев\cite{sob54}, \,О.А.Ладыженская\cite{lad},  К.Фридрихс \cite{fri},
 	Э.Быховский и Н.Смирнов \cite{bs}
 	приводят аналогичные разложения $\mathbf{L}_{2}(G)$ на  
 	ортогональные подпространства.  Так,  С.Л.Соболев 
 	предполагает, что область
 	$\Omega$ гомеоморфна шару. В этом случае  пространство $\mathcal{B}_H(\Omega)$  пусто.
 	
 	О.А.Ладыженская  
 	приводит  разложение: $\mathbf{L}_{2}(\Omega)=\mathbf{G}(\Omega)\oplus \mathbf{J}^\circ(\Omega) $, 
 	где   $\mathbf{J}^\circ (\Omega) $ замыкание  в норме  $\mathbf{L}_{2}(\Omega)$  множества бесконечно дифференцицируемых финитных  в  $\Omega$ соленоидальных векторов, а  $\mathbf{G}(\Omega)$  состоит из $grad\, \varphi $,   где  $\varphi$ есть  однозначная в   $\Omega$  функция, локально квадратично суммируемая и имеющая первые производные из  $\mathbf{L}_{2}(\Omega)$ (Теорема 1 $\S 2$ Гл. 1 в \cite{lad}).     } 
 Мы   будем придерживаться разложения (1.1). 

  \subsection{  Операторы градиент, ротор и дивергенция} Определяются в трехмерном
векторном анализе \cite{zo}. Им соответствует оператор $d$ внешнего
дифференцирования на формах $\omega^k$ степени $k=0,1$ и 2.
Соотношения $dd\omega^k=0$ при $k=0,1$ имеют вид $\mathrm{rot}\,\nabla h=0$ и
 $\mathrm{div}\, \mathrm{rot} \mathbf{u}=0$.\quad 
Формулы $\mathbf{u}\cdot\nabla
h+h\mathrm{div}\mathbf{u}=\mathrm{div}(h \mathbf{u}) $, \quad 
$\mathbf{u}\cdot\mathrm{rot} \mathbf{v}- \mathrm{rot}
\mathbf{u}\cdot\mathbf{v}=\mathrm{div}[\mathbf{v},\mathbf{u}]$, где
$[\mathbf{v},\mathbf{u}]$ - векторное произведение, и интегрирование
по области $G$ используются при определении операторов $\mathrm{div}
\mathbf{u}$ и $\mathbf{rot}\mathbf{u}$ в $\mathbf{L}_{2}(G)$.
 Оператор Лапласа
выражается через $\mathrm{rot}\, \mathrm{rot}$ и $\nabla \mathrm{div}$:
\begin{equation*}\klabel{dd 1}\Delta \mathbf{v} =\nabla \mathrm{div} \mathbf{v}
-\mathrm{rot}\, \mathrm{rot} \mathbf{v}. \eqno(1.3) \end{equation*}
 Оператор Лапласа эллиптичен  \cite{ so71},  а операторы $\mathrm{rot}$ и
 $\nabla \mathrm{div}$ не являются таковыми. Они
вырождены, причем $\mathrm{rot}\, \mathbf {u}=0$ при $\mathbf{u}\in~
\mathcal{A}$, $\nabla\mathrm div \mathbf {v}=0$ при $\mathbf{v}\in
\mathcal{B}$ в смысле  $\mathbf{L}_{2}(G)$ \cite{hw}.

Поэтому $\Delta \mathbf{v} = \nabla \mathrm{div} \mathbf{v}$ на
$\mathcal{A}$ и $\Delta \mathbf{v} =-\mathrm{rot}\, \mathrm{rot}
\mathbf{v} $ на $\mathcal{B}$.

{\it Замечание.}{ \small H.Weyl называет безвихревым (irrotational)  поле  $\mathbf{u}\in \mathbf{L}_{2}(G)$,  для которого  $(\mathbf{u}, \mathrm{rot}\, \mathbf{v}) =0$ для  любого поля	$ \mathbf{v}$ c компонентами $v_j$  из $\mathcal{C}_0^1(G)$,  а поле $\mathbf{w}\in \mathbf{L}_{2}(G)$,  для которого  $(\mathbf{w}, \nabla\, v) =0$ 		 $\forall \,v\in \mathcal{C}_0^1(G)$, - соленоидальным \cite{hw}. \quad  Символ  "$\mathrm{rot}\, \mathbf {u}=0$ при $\mathbf{u}\in\mathcal{A}"$   означает,  что $\mathbf{u}=\{\nabla h\}$,  где $ h\in H^1(G)$,  и  $(\mathbf{u}, \mathrm{rot}\, \mathbf{v})=  (\nabla h, \mathrm{rot}\, \mathbf{v})=0$ для  любого	$ \mathbf{v}$ c компонентами $v_j\in \mathcal{C}_0^{\infty}(G)$,  что очевидно.  }
 \subsection{ Содержание: классы обобщенно элиптических задач } В $\S 1$ определяются основные подпространства  $\mathbf{L}_{2}(G)$,  операторы, их соотношения,  и формулируются основные результаты работы; 
 $\S 2$  содержит  постановку   краевых   задач  (2.1), (2.2) для операторов
 $ \text{rot}+\lambda \, I$    и  $ \nabla\text{div}+\lambda \,I$     первого и второго порядков в  пространствах Соболева.  
  Определяются  классы [REESp] обобщенно эллиптических систем. 
Системы (2.1), (2.2) при $\lambda\neq 0$  принадлежат  классу  [REES1].
 Им  соответствуют операторы $\mathbb{A}$ и $\mathbb{B}$.  
   
   Мы доказываем, что задачи (2.1) и (2.2) удовлетворяют 
 условиям эллиптичности В.Солонникова и применяем его Теорему 1.1 \cite{so71}.   Из нее вытекают
 \begin{theorem} \klabel{rot   1}
	Оператор $\mathbb{A}$	имеет левый регуляризатор.  
	Его ядро конечномерно и для любых $ \mathbf{u}\in  \mathbf{H}^{s+1}(G)$ и $ \lambda \neq 0 $ ( с  постоянной  $C_s =C_s(\lambda )>0$,  зависящей только от $s, \lambda $)	выполняется оценка: 
 \begin{equation*}
\klabel{arot__1_} C_s\|\mathbf{u}\|_{s+1}
\leq\|\mathrm{rot} \mathbf{u}\|_{s}+
|\lambda|\|\mathrm{div} \mathbf{u}\|_{s}+
|\gamma({\mathbf{n}}\cdot\mathbf{u})|_{s+1/2}+
\|\mathbf{u}\|_{s}.\eqno{(1.4)}
\end{equation*}
\end{theorem} 
  \begin{theorem} \klabel{Nd_1}
  	Оператор $\mathbb{B}$ 	имеет левый регуляризатор.
  Его ядро конечномерно и   для любых $\mathbf{u}\in  \mathbf{H}^{s+2}(G)$ и $\lambda \neq 0 $ ( с  постоянной  $C_s =C_s(\lambda )>0$,  зависящей только от $s, \lambda $)	выполняется оценка:
 \begin{equation*}
\klabel{ond__2_} C_s\|\mathbf{v}\|_{s+2}
\leq|\lambda|\|\mathrm{rot} \mathbf{v}\|_{s+1}+
\|\nabla\mathrm{div} \mathbf{v}\|_{s}+
|\gamma({\mathbf{n}}\cdot\mathbf{v})|_{s+3/2}+ \|\mathbf{v}\|_{s}.\eqno{(1.5)}
\end{equation*}
\end{theorem}
 Существание левого регуляризатора  оператора  $\mathbb{B}$ и оценка (1.5)  доказаны  в \cite{ds18}.

  Топологических ограничений на область нет, предполагается ее
 связность, ограниченность и гладкость границы.
 
 Оценка (1.5)  получена мною {\it впервые} из оценок Солонникова [14].
 Оценка (1.4) известна  мне давно,  она не была выписана  в \cite{sa72},  хотя я знал, что  элиптичность задачи эквивалентна точной оценке в пространствах Соболева от Л.Р.Волевича \cite{vol},  тогда я работал еще в пространствах   Гельдера.  
  Joshida и Giga в  [12] ссылаются на работы J.P.Bourguignon, H. Brezis \cite{bobr} и C. Foias, R.Temam \cite{foit}

 Этот подход   применим для других обобщенно эллиптических систем класса [REESp]. Этот класс  выделен из класса Вайнберга и Грушина  \cite{vagr},  он содержит   системы,  главные части  которых суть степени ротора или градиента дивергенции. 

Из эллиптической теории вытекают   свойства решений спектральных задач операторов    ротора и градиента дивергенции,  такие как  конечная кратность ненулевых с.-значений и   гладкость с.-полей в  любой области $G$ с гладкой границей. 

   Явные    формулы решений спекральных задач  операторов    
   ротора и градиента дивергенции в шаре    \cite{saUMJ13} я нашел,  обнаружив их связи с решениями   спекральных задач Дирихле и Неймана для оператора Лапласа.  Они решены в учебнике   В.С. Владимирова \cite{vla}.
 \subsection{ Содержание: класс  $\mathbf{V}^{0}$  вихревых полей} 
 Собственные поля ротора образуют базис в   $\mathbf{V}^{0}(G)$.  В $\S 3$ определяются также пространства   
$\mathbf{W}^{k}(G)$, - (3,7), и операторы  $S$, - (3.5),  $S^{-1}$, -  (3.9):    
\[\mathbf{W}^{k}(G)=
\{\mathbf{f}\in\mathbf{V}^{0}(G),..., \mathrm{rot}^k
\mathbf{f}\in\mathbf{V}^{0}(G) \}\quad \text{при}\quad k\geq
1.\]
  Оператор
$S^{-1}$ отображает   $\mathbf{V}^{0}$  на  $\mathbf{W}^{1}$, а  $\mathbf{W}^{k-1}$ на $\mathbf{W}^{k}$  при  $k>1$.
Оператор
$S$ отображает  $\mathbf{W}^{k}$ на $\mathbf{W}^{k-1}$, 
а  $\mathbf{W}^{1}$ на $\mathbf{W}^{0}\equiv \mathbf{V}^{0} $.

Рассматривается также оператор $S+\lambda I$.  Мы доказываем,  что

 Оператор $S+\lambda I:\mathbf{W}^{k}	\rightarrow \mathbf{W}^{k-1}$ - фредгольмов.
По определению оператор $ S+\lambda I$ совпадает с 
$\mathrm{rot}+\lambda I$   на $\mathbf{W}^{1}$  и 
\[(S+\lambda I)\mathbf{f}=\lim_{n\rightarrow\infty}
(\mathrm{rot+\lambda I}) (\mathbf{f}^n_{\mathbf{V}})
=\sum_{j=1}^\infty[(\lambda+\lambda_j)(\mathbf{f},\mathbf{q}^{+}_{j})\mathbf{q}^{+}_{j}+
(\lambda-\lambda_j)(\mathbf{f},\mathbf{q}^{-}_{j})\mathbf{q}^{-}_{j}].\]
Ряд сходится в $\mathbf{L}_{2}(G)$,  так как $\| (S+\lambda I)\mathbf{f}\|^2_{ \mathbf{V}^0}\leq c^2_0 \|\mathbf{f}\|^2_{\mathbf{W}^1}$, где $c_0<\infty  $ (см. (3.14)).
Обратный оператор имеет вид:
\begin{equation*}\klabel{sp__3_}
(S+\lambda I)^{-1}\mathbf{f}=
\sum_{j=1}^{\infty}\left[
\frac{(\mathbf{f},{\mathbf{q}}_{j}^+)}{\lambda+\lambda_{j}}
\mathbf{q}_{j}^+(\mathbf{x})+
\frac{(\mathbf{f},{\mathbf{q}}_{j}^-)}{\lambda-\lambda_{j}}
\mathbf{q}_{j}^-(\mathbf{x})\right], \eqno (1.6) \end{equation*}
если ни одно из слагаемых этого ряда не обращается в
бесконечность.  Это означает, что либо $\lambda\pm \lambda_{j}\neq 0$ для всех $j$, либо  $(\mathbf{f},\mathbf{q}^{-}_{j})=0$ при  $\lambda=\lambda_j=\lambda_{j_0}$  и эти элементы отсуствуют в ряду.\quad
 При этом  \newline
$ \|(S+\lambda I)^{-1}\mathbf{f}\|^2_{\mathbf{W}^1} \leq C^2_0 \|\mathbf{f}\|^2_{\mathbf{V}^0}$,  где 
$C^2_0 <\infty$ не зависит от $\mathbf{f}$   (см. (3.16)).
Следовательно, оба оператора непрерывны и имеет место
\begin{theorem}\klabel{S  1}
Оператор $S+\lambda I:
	\mathbf{W}^{1}(G)\rightarrow \mathbf{V}^{0}(G)$ непрерывен и
	однозначно обратим, если $\lambda$ не принадлежит спектру
	$\sigma_p(S)\subset\mathbb{R}$
	оператора $S$. Его обратный задается формулой
	(1.6) и для любого  $\mathbf{f}\in {\mathbf{V}^0}$  ряд $(S+\lambda I)^{-1}\mathbf{f}\in {\mathbf{W}^1}$ .\quad 	
	
	Если $\lambda=\lambda_{j_0}$, то он обратим тогда и только
	тогда, когда
	\begin{equation*}
	\klabel{urz _1_}\int_G \mathbf{f}\cdot \mathbf{q_j^-} dx=0\quad
	\text{для}\quad\forall \mathbf{q_j^-}: \lambda_j=\lambda_{j_0}. \eqno{(1.7)}
	\end{equation*}
	Ядро оператора $S+\lambda_{j_0} I$ конечномерно и определяется
	собственными функциями $\mathbf{q_j^-}(\mathbf{x})$, собственные
	значения которых равны $\lambda_{j_0}$:
	\begin{equation*} \klabel{ker__1_}
	Ker(S+\lambda_{j_0} I)= \sum_{\lambda_j=\lambda_{j_0}}
	c_j \mathbf{q}^{-}_{j}(\mathbf{x} \quad\forall \,c_j\in
	\mathbb{R}.\eqno{(1.8)} \end{equation*}\end{theorem}
В п. 3.5 приводятся  также  оценки  
\[ \|  (S+\lambda I)\mathbf{f}\|^2_{ \mathbf{W}^m}\leq  c^2_m \|\mathbf{f}\|^2_{\mathbf{W}^{m+1}}, \,\,  \|(S+\lambda I)^{-1}\mathbf{f}\|^2_{\mathbf{W}^{m+1}}\leq  C^2_m \|\mathbf{f}\|^2_{\mathbf{W}^{m}},  \eqno{(1.9)} \]
где  постоянные $  c_m,\,  C_m< \infty $ не зависят от $\mathbf{f}$ и  при $m=0$ совпадают с (3.14) и (3.16).
 Из теоремы и этих оценок следует 
\begin{lemma}\klabel{S}  Если $\lambda \overline{\in} Sp(\mathrm{S})$, 	 $ m \geq0$, то 	 операторы  $\mathrm{S}+\lambda I$ ( и его обратный)  отображают    пространство  $ \mathbf{W}^{m+1}$	на  $  \mathbf{W}^{m}$	 (и обратно) взаимно	однозначно и непрерывно. \end{lemma}
 \subsection{ Соотношения между пространствами  $\mathbf{W}^{k}$,
$\mathbf{H}^{k}$ и $\mathbf{C}^{k-2}$. } Область
$\Omega$, гомеоморфную шару, выделил С.Л.Соболев в [4]. В этом случае  пространство $\mathcal{B}_H(\Omega)$  пусто. Границу  области $\Omega$  будем предполагать гладкой.
Скалярное произведение
в $\mathbf{H}^k(\Omega)$ Сергей  Львович определяет так:
\begin{equation*} \klabel{skpro 1}
(\mathbf{f},\mathbf{g})_k=(\mathbf{f},\mathbf{g})+
\int_{\Omega} \sum_{|\alpha|=k}\frac{k!}{\alpha
	!}\partial^{\alpha}\mathbf{f}\cdot\partial^{\alpha}\mathbf{g}
d \mathbf{x},\quad k\geq 1. \quad \eqno{(1.10)}
\end{equation*} В  пространстве $ \mathbf{W}^k(\Omega)=
\{\mathbf{f}\in \mathbf{V}^0,
..., \mathrm{rot}^{k} \mathbf{f}\in \mathbf{V}^0\}$
норму $\mathbf{f}\in{\mathbf{W}^k}$ выбираем   так же:
$\|\mathbf{f}\|^2_{\mathbf{W}^k}\equiv \|\mathbf{f}\|^2+\|\mathrm{rot}^k \mathbf{f}\|^2$.
\quad Имеет место
\begin{theorem}\klabel{W}
	Для того, чтобы $\mathbf{f}\in \mathbf{V}^0(\Omega)$ разлагалась в ряд Фурье
	\begin{equation*} \klabel{rof 1}
	\mathbf{f}(\mathbf{x})=\sum_{j=1}^{\infty}
	((\mathbf{f},\mathbf{q}_{j}^+)\mathbf{q}_{j}^+(\mathbf{x})
	+(\mathbf{f},\mathbf{q}_{j}^-)\mathbf{q}_{j}^-(\mathbf{x})), \quad
	\|\mathbf{q}_{j}^{\pm}\| =1,\eqno{(1.11)}
	\end{equation*}
	по собственным вектор-функциям $\mathbf{q}_{j}^{\pm}(\mathbf{x})$
	ротора в области $\Omega$,
	сходящийся в норме
	пространства Соболева $\mathbf{H}^k(\Omega)$, необходимо и достаточно,
	чтобы $\mathbf{f}$ принадлежала $\mathbf{W}^k(\Omega)$.
	
	Если $\mathbf{f}\in \mathbf{W}^k(\Omega)$,
	то существует такая постоянная $C>0$, не зависящая от
	$\mathbf{f}$, что
	\begin{equation*} \klabel{orf 3}
	\sum_{j} {\lambda}_{j}^{2k} ((\mathbf{f},\mathbf{q}_{j}^+)^2
	+(\mathbf{f},\mathbf{q}_{j}^-)^2)\leq C
	\|\mathbf{f}\|^2_{\mathbf{H}^k(\Omega)}.\eqno{(1.12)}
	\end{equation*}
		Если $k\geq 2$, то  вектор-функция $\mathbf{f}$ из	$\mathbf{W}^k(\Omega)$
	разлагается в в ряд (1.11), сходящийся в пространстве 
	$\mathbf{C}^{k-2}(\overline{\Omega})$.     \end{theorem}
{\bf Следствие.} {\it Вектор-функция $\mathbf{f}\in
	\mathbf{V}^0\cap\mathbf{C}^{\infty}_0(\Omega)$
	разлагается в	ряд	(1.11),
	сходящийся в любом из пространств
	$\mathbf{C}^{k}(\overline{\Omega})$, \, $ k\in \mathbb{N}$.}

Эти результаты дополняют известные в теории рядов Фурье утверж-дения
(см. Теорема 7 в $\S 4$ гл. 2 в [5], Теорема 8 в $\S 2$ гл. 4 в [2]).
Таким образом, $\mathbf{W}^k(G)$ - это  аналоги
пространств Соболева $\mathbf{H}^k(G)$ в классе 
соленоидальных  полей.\quad  
Отметим вложения: \[\mathbf{W}^k\subset \mathbf{W}^{k-1}\subset...\subset \mathbf{W}^1\subset \mathbf{V}^{0}.  .\eqno{(1.13)} \]

Заметим, что   пространств $\mathbf{W}^{k}$,   $k> 1$,   Joshida и Giga   \cite{yogi} не рассматривали.  Они ввели ${H}^{1}_{\Sigma\Sigma}=\mathcal{D}(S)$ как область определения  $S$.  
  Пространства $ \mathbf{W}^k$ 
   и     соотношения   между  ними и
    $\mathbf{H}^{k}$ и $\mathbf{C}^{k-2}$ (Теоремы 3,  4 и Лемма 1) -   это {\bf первый основной} результат  этой статьи.
  
 \subsection{ Класс  $\mathcal{A}$ потенциальных  полей }
 Недавно   (8.02) вышел в свет очередной номер журнала  Динамические системы,  2018, т. 8, где в статье  \cite{ds18}
 "Оператор градиент дивергенции  и  пространства Соболева"  \,  изучен
  класс  $\mathcal{A}$ потенциальных  полей: собственные поля оператора $\nabla\mathrm{div}$ задают ортогональный базис в  $\mathcal{A}_{\gamma}$,   оператор  $\mathcal{N}_d$  есть  самосопряженное расширение $\nabla\mathrm{div}$  в   $\mathcal{A}_{\gamma}$,    пространства 
\[ \mathbf{A}_{\gamma}^{2k}(G)=
\{\mathbf{f}\in\mathcal{A}_{\gamma}(G),...,
(\nabla\mathrm{div})^k\,\mathbf{f}\in\mathcal{A}_{\gamma}(G) \}\quad
\forall \, k\geq 1,    \eqno{(1.14)}\]
-- аналоги пространств Соболева $\mathbf{H}^{2k}(G)$
порядков  $2k$  в классе  $\mathcal{A}_{\gamma}$.

 Так же как Лемма 1 доказана
\begin{lemma}\klabel{N} 
Если $\lambda \overline{\in} Sp (\mathcal{N}_d)$,  $k \geq0$,
то 	 операторы  $ \mathcal{N}_d+\lambda I$  (и его обратный)  отображает    пространство  $ \mathbf{A}^{2(k+1)}_{\gamma }$
на  $  \mathbf{A}^{2k}_{\gamma }$ (и обратно) 	взаимно	однозначно и непрерывно.
\end{lemma}
Отметим вложения \[\mathbf{A}^{2k}_{\gamma }\subset ...\subset\mathbf{A}^{2}_{\gamma }\subset\mathcal{A}_{\gamma }\subset \mathcal{A}\subset \mathbf{L}_2(G). \eqno{(1.15)}\]

 Базисные векторы  в классах  $\mathcal{A}$ и $\mathcal{B}=\mathcal{B}_H\oplus \mathbf{V}^0$  в совокупности образуют базис во всем пространстве  $\mathbf{L}_{2}(G)$.

 \subsection{Содержание: классы пространств $ C(2k,m)$   в $\mathbf{L}_{2}(\Omega)$ }  
В $\S 5$ рассматривается 
 целочисленная решетка  
пространств $ C(2k,m)\equiv \mathbf{A}^{2k}_{\gamma } \oplus \mathbf{W}^m$,  называемых классами, $k\geq 0, m\geq 0 $  -целые, $k+m>0$ .  А также пространство 
	$\mathbf{E}_{\gamma}^{0}(\Omega)$.
Доказаны

\begin{theorem}\klabel{R}
	  Если $\lambda\neq 0, \pm {\lambda}_{j}$,
	$j\in \mathbf{N}$ и
	$\mathbf{f}\in\mathbf{L}_{2}(\Omega)$,
	то
	единственное решение $\mathbf{u}$ задачи 1 дается суммой рядов-проекций
	$\mathbf{u}_{\mathcal{A}}+\mathbf{u}_\mathbf{V}$:
	\begin{equation*} \klabel{kra 1}
	{\mathbf{u}_{\mathcal{A}}}={\lambda}^{-1}\mathbf{f}_\mathcal{A}\equiv{\lambda}^{-1}
	\sum_{j=1}^{\infty} (\mathbf{f},{\mathbf{q}}_{j}) \mathbf{q}_{j}
	(\mathbf{x}),\eqno{(1.16)}
	\end{equation*}
	\begin{equation*} \klabel{kro 2}
	\mathbf{u}_\mathbf{V}=(S+\lambda)^{-1} \mathbf{f}_\mathbf{V}\equiv
	\sum_{j=1}^{\infty}\left[\frac{(\mathbf{f},{\mathbf{q}}_{j}^{+})}
	{\lambda{+}\lambda_{j}} \mathbf{q}_{j}^{+}(\mathbf{x})+
	\frac{(\mathbf{f},{\mathbf{q}}_{j}^{-})}{\lambda-\lambda_{j}}
	\mathbf{q}_{j}^{-}(\mathbf{x})\right]. \eqno{(1.17)}\end{equation*}
	
	В частности, если
 $\mathbf{f}=\mathbf{f}_{\mathcal{A}}$ и $\mathbf{f}_{\mathcal{A}}\in \mathcal{A}$ или  $\mathbf{f}_{\mathcal{A}}\in\mathcal{A}_{\gamma}$,
 
то $\mathbf{u}=~{\lambda}^{-1}\mathbf{f}_{\mathcal{A}}\in
\mathcal{A}$  или   $\mathbf{u}\in \mathcal{A}_{\gamma} $ -  обобщенные решения  задачи 1.

	Если $\mathbf{f}\in \mathcal{B}\bot \mathcal{A}$ в
	$\mathbf{L}_{2}(\Omega)$, то
	$\mathbf{u}=(S+\lambda)^{-1}\mathbf{f}_\mathbf{V}\in \mathbf{W}^1\subset \mathbf{H}^{1}_{\gamma}(\Omega)$.
	
	Если 	$\mathbf{f}\in\mathbf{E}_{\gamma}^{0}(\Omega)$,  то 	$\mathbf{u}=\mathbf{u}_{\mathcal{A}}+\mathbf{u}_\mathbf{V}\in  \mathbf{H}^{1}_{\gamma}(\Omega)$.
	
	Если 	$\mathbf{f}$  принадлежит классу $C(2k, m)$,  то	$\mathbf{u}\in C(2k, m+1)$.

	Если	же $\mathbf{f}\in \mathcal{D}(\Omega)$, то ряды (1.16),
	(1.17) сходятся в 	$\mathbf{H}^{s}(\Omega)$
	
	 для 	любого $s\geq 1$ и $\mathbf{u}\in
	C^{\infty}(\overline{\Omega})$ - классическое решение задачи.
\end{theorem}
{\it Замечание.}
{\small  В статье   \cite{saUMJ13} мы доказали,  что собственные значения   ротора  в шаре радиуса $R$ равны $\pm \rho_{n,m}/R$,  где числа $\pm \rho_{n,m}$ - нули функций 
		\begin{equation*}\klabel{psi__1_}
	\psi_n(z)=(-z)^n\left(\frac{d}{zdz}\right)^n\frac{\sin
		z}z,\quad m,\, n\in \mathbb {N}, \eqno{(1.18)}\end{equation*}	
	 кратность  собственного значения  $\lambda_{n,m} =\pm \rho_{n,m}/R$ равна $2n+1$.	
	 
	Собственные значения  оператора $\nabla\mathrm{div}$ равны $-\nu_{n,m}^2$, где $\nu_{n,m}=\alpha_{n,m}/R$,  а  числа  $\alpha_{n,m}$ - нули производных  $\psi'_n(r)$,   $n \geq 0,\, m\in {\mathbb {N}}$ ;
	
	 кратность  собственного значения  $-\nu^2_{n,m}$ равна $2n+1$.
	
	Собственные поля $\mathbf{q}_{\kappa}$  градиента дивергенции и $ \mathbf{q}^{\pm}_{\kappa}$
	ротора  выражаются явно через сферические функции и  функции $\psi_n(r)$;  $ \kappa=(n,m,k)$.}

Из Теоремы 5 и Леммы 1 вытекают 
\begin{lemma}\klabel{Ro}
	  При $\lambda \overline{\in} Sp (\mathrm{rot})$ оператор  $\mathrm{rot}+\lambda I$ 
	отображает класс  $C(2k, m+1)$
	на  $C(2k, m)$	взаимно	однозначно и непрерывно, $k,m\geq~0$. \end{lemma}
	
	{\bf Следствие. } {\it
 Если область $\Omega=B$,  
	  $\psi_n(\lambda\,R)\neq 0$   $\forall\, n\in {\mathbb {N}}$, а поле  $\mathbf{f}\in  \mathbf{A}^{2k}_{\gamma }(B)\oplus  \mathbf{W}^m(B)$, то решение задачи 1 существует, единственно и принадлежит классу    $ \mathbf{A}^{2k}_{\gamma }(B)\oplus  \mathbf{W}^{m+1}(B)$.  }

Аналогично доказаны 
	\begin{lemma}\klabel{o}   При     $\nu^2\neq Sp(-\nabla \, \mathrm{div})$ оператор  $ \nabla \,  \mathrm{div}+\nu^2\,I$  отображает  $C(2(k+1), m)$ на класс $C(2k,m)$ взаимно однозначно и непрерывно. \end{lemma}
   	{\bf Следствие. } {\it Если область $\Omega=B$,
  $\psi'_n(\nu\,R)\neq 0$ $\forall\, n \geq 0$, а поле  $\mathbf{f}\in  \mathbf{A}^{2k}_{\gamma }(B)\oplus  \mathbf{W}^m(B)$, то решение задачи 2 существует, единственно и принадлежит классу    $ \mathbf{A}^{2(k+1)}_{\gamma }(B)\oplus  \mathbf{W}^{m}(B)$. }

Таким образом, изучены   пространства  $\mathbf{W}^m$ на рядах Фурье, определяемых собственными полями оператора ротор (вихрь).
В пространстве  $ \mathbf{L}_{2}(\Omega)$ введены классы $ C(2k,m)\equiv  \mathbf{A}^{2k}_{\gamma } \oplus \mathbf{W}^m$   и рассмотрены их отображения операторами  $\mathrm{rot}+\lambda I$ и   $ \nabla \, \mathrm{div}+\nu^2\,I$.\newline 
Теорема 5,   Леммы 3, 4 и их следствия - это {\bf второй основной} результат автора этой статьи.

\section { Краевые и спектральные задачи }
 \subsection { Краевые задачи}
 
 В ограниченной области $G$ с гладкой границей
$\Gamma$  изучаются {\it задачи:} найти вектор-функции
$\mathbf{u}$ и $\mathbf{v}$, такие что
$$
\klabel{rot__1_}
 \text{rot}\mathbf{u}+\lambda \mathbf{u}=\mathbf{f}(\mathbf{x}), \quad
 \mathbf{x}\in G, \quad \mathbf{n}\cdot \mathbf{u}{{|}_{\Gamma }}=g,
 \eqno{(2.1)}
$$
\begin{equation*}
\klabel{grd__1_}
 \nabla\text{div}\mathbf{v}+\lambda \mathbf{v}=\mathbf{f}(\mathbf{x}), \quad
 \mathbf{x}\in G, \quad \mathbf{n}\cdot \mathbf{v}{{|}_{\Gamma }}=g,
 \eqno{(2.2)}
\end{equation*} где векторная и
скалярная функции $\mathbf{f}$ и ${g}$ заданы. Решения задач ищем в
пространствах Соболева $\mathbf{H}^{s+1}(G)$ и
$\mathbf{H}^{s+2}(G)$, где $s$-целое $\geq 0$, а $(\mathbf{f},{g})$
задаем в пространствах: $\mathbf{f}\in \mathbf{H}^{s}(G), g\in
{H}^{s+1/2}(\Gamma)$ и $\mathbf{f}\in \mathbf{H}^{s}(G), g\in
{H}^{s+3/2}(\Gamma)$, соответственно.
Эта постановка является {\it классической} в теории {\it
эллиптических краевых задач в пространствах Соболева:} см.
 Соболев [1], Солонников [14].

Отметим, что ненулевые решения $(\mathbf{u}, \lambda)$ и
$(\mathbf{v}, \lambda)$ однородных задач (2.1) и (2.2) ( $\mathbf{f}=0$ и ${g}=0$)- это решения
спектральных задач операторов $ \text{rot}$ и $\nabla\text{div}$.\quad 
Они аннулируют друг друга и
\begin{equation*}\klabel{anrot__1_}
\mathrm{rot}\mathbf{u}=0 \,\text{на} \, \mathcal{A}=\{\nabla h,
 \,h \in H^{1}\}, \quad \nabla \mathrm{div}\mathbf{v}=0 \quad
 \text{на} \quad \mathcal{B}\perp \mathcal{A}. \end{equation*}
 Ортогональные пространства
${\mathcal{{A}}}$ и ${\mathcal{B}}$ в $\mathbf{L}_{2}(G)$
 бесконечномерны [3].
При $\lambda=0$ однородные задачи (2.1) и (2.2) имеют счетное число линейно независимых решений.
Значит, {\it нулевая точка} спектра каждого из операторов $
\text{rot}$ и $\nabla\text{div}$ имеет {\it бесконечную кратность}.
Специфика этих задач состоит в том, что 
эти операторы при $\lambda\neq 0$ являются {\bf обобщенно эллиптическими} класса  [REES 1]. 
 \subsection{ Класс   систем, приводимых к эллиптическим системам}
Определение этого класса  мы приведем  для
систем дифференциальных уравнений с постоянными коэффициентами,

Система дифференциальных уравнений, $S(D)u=f$ порядка $m$, из этого
класса обладает свойствами: \newline
а) ее символическая матрица
$S_0(i\xi)$ имеет постоянный ранг для любой  $\xi\in\mathcal{R}^3\backslash 0$. \quad 
Это позволяет построить
аннулятор $C(D)$ оператора $S_0(D)$ такой, что $(CS_0)(D)\equiv 0$
на $X$ и определить
\newline б) расширенную систему $\left(\begin{matrix}Su=f\\ CSu=Cf\end{matrix}\right)$
  порядка   $\left(\begin{matrix}m\\k\end{matrix}\right)$.
\newline
Ее символическая матрица
$\left(\begin{matrix}S_0(i\xi)\\(CS)_0(i\xi)\end{matrix}\right)$
определяется
младшей частью оператора $S(D)$ и дополняет матрицу $S_0(i\xi)$.
\newline в) Если ранг расширенной матрицы максимален, то
исходная система $Su=f$ принадлежит классу [REES 1] и степень ее
приводимости равна единице.\newline 
г) Если система $Su=f$ такова, что ранг расширенной матрицы не
максимален, но постоянный, то процесс повторяется и при определенных
условиях система принадлежит классу [REES 2].... Символ 
[REES p] означает: REduced to Elliptic Systems на  p -том шаге. 
   
    Б.Вайнберг и В.Грушин [15] доказали, что система $Su=f$ класса [REES p] является
разрешимой по Фредгольму или Нетеру в пространствах Соболева
$\mathbf{H}^s(X)$, если $f\in \mathbf{H}^{s-m+p}(X)$, где $s\geq m$ целое.
 В качестве примера 
  приводится оператор $d+\ast$ на дифференциальных формах степени $k$ в
$2k+1$-мерном многообразии $X$ без края.

Системы  (2.1) и (2.2) являются {\bf обобщенно эллиптическими} класса  [REES 1].
  Действительно, из этих уравнений  вытекают соотношения:  $
\lambda \text{div}\,\mathbf{u}=\text{div}\,\mathbf{f}(\mathbf{x}), \quad
\lambda \text{rot}\,\mathbf{v}=\text{rot}\,\mathbf{f}(\mathbf{x}),
\quad \mathbf{x}\in G.$
 Соединяя их в систему, видим, что операторы
\begin{equation*}
\klabel{ell__12_}\left(\begin{matrix}\text{rot}+\lambda I \\
\lambda \text{div}\end{matrix}\right)\quad \text{и}\quad
\left( \begin{matrix} \nabla\text{div}+\lambda I\\
\lambda \text{rot}\end{matrix}\right)\eqno{(2.3)}
\end{equation*} являются эллиптическими по Даглису-Ниренбергу [14].
Значит, они принадлежат классу  [REES 1] систем дифференциальных уравнений,
 приводимых к эллиптическим системам на  первом шаге расширений  Вайнберга
и Грушина [15].  
\footnote {Другие классы обобщенно
	эллиптических операторов см. в работе \cite{saks97}.}
 \subsection{  Обобщенно эллиптическая краевая задача.} Рассмотрим подробнее первую из них.
 Расширенная система
\begin{equation*} \label{rodi__1_}
\mathrm{rot}\,\mathbf{u}+\lambda \mathbf{u}=\mathbf{f},\quad \lambda
\mathrm{ div}\,\mathbf{ u}=\mathrm{div}\,\mathbf{f}, \eqno{(2.4)}\end{equation*}
является эллиптической системой первого порядка
( переопределеной, если $f_4\neq \mathrm{div}\,\mathbf{f}$).
Вместе с краевым условием
$\gamma\,\mathbf{n}\cdot \mathbf{u}=g$
она составляют    эллиптическую краевую задачу
по Солонникову \cite{so71}. Это означает, что \newline 
1) система  (2.4) эллиптична,\footnote {  Главные части  системы в  \cite{so71}
	определяются с помощью весов $s_k$ и $t_j$ таких, что
	$ord\, L_{k,j}\leq s_k+t_j $. 
	Положив    $s_k=0$\, при $k=1,...,4$
	и $t_j~=1$ при $j=1,2,3$ мы получим операторы
	системы (2.4),  а в краевом операторе положим	$\sigma_1=-1$.}
\newline 2)  краевое условие $\gamma\mathbf{n}\cdot \mathbf{u}$\,
"накрывает" \,\,оператор системы. 

Первое условие сводится к тому, что  однородная система линейных
алгебраических уравнений:
\begin{equation*}
\label{cdx  1}\mathrm{rot}(i\xi )\mathbf{w}=0, \quad
\lambda\,\mathrm{div}(i\xi )\mathbf{w}=0, \quad \forall \xi\neq 0 \eqno{(2.5)}
\end{equation*}
c параметром $\xi \in \mathcal{R}^3$ имеет только тривиальное
решение $ \mathbf{w}=0$.

Второе условие означает, что  однородная система линейных
диф-ференциальных уравнений:
\begin{equation*}\label{cdz  2}
\mathrm{rot}(i\tau+\mathbf{n} d/dz ) \mathbf{v}=0,\quad
\mathrm{div}(i\tau+\mathbf{n} d/dz )\mathbf{v}=0, \quad \forall \tau
\neq 0,     \eqno{(2.6)} \end{equation*} на полуоси $z\geq 0$ с краевым условием:
$\mathbf{n}\cdot \mathbf{v}|_{ z=0}=0$ и убыванием,
$\mathbf{v}(y,\tau; z )\rightarrow 0$ при $z\rightarrow + \infty$,
имеет только тривиальное решение.

Здесь  $\tau$ и $\mathbf{n}$ -- касательный и нормальный векторы к
$\Gamma$ в точке $y\in \Gamma$ и $|\mathbf{n}|=1$.
При доказательстве этих утверждений   воспользуемся соотношением 
\quad $ \mathrm{rot}\, \mathrm{rot}\,\mathbf{v}= -\Delta \mathbf{v}
+ \nabla \mathrm{div} \mathbf{v}.$

Тогда $1^0)$. Из уравнений $ (2.5)$ вытекает уравнение
$-\Delta(i\xi)\mathbf{w}=0$. Оно распадается на три скалярных
уравнения $|\xi | ^2 w_j =0$. 
Значит, $\mathbf{w}=0$ при $|\xi | \neq 0$. Эллиптичность системы
(2.4)  доказана.
\newline
$2^0)$. Из уравнений $(2.6)$ получаем уравнение $ (-|\tau| ^2
+ (d/dz)^2)\mathbf{v} = 0$ с параметром $|\tau |> 0$. Его убывающее
при $z\rightarrow + \infty$ решение имеет вид:
$\mathbf{v}=\mathbf{w} e^{-|\tau|z}$ . Оно удовлетворяет уравнениям
$(2.6)$, если вектор-функция $\mathbf{w}$ такова, что
$ \omega\times \mathbf{w}=0,\quad {\omega}' \cdot
\mathbf{w}=0$,   где $ \omega\equiv i\tau -| \tau| \mathbf{n}$
--вектор-столбец,  $ \omega'$ -- вектор-строка, а ${\omega}' \cdot
{\omega}$ -- их произведение.

Легко убедиться, что векторное и  "скалярное"  произведения $\omega$
на $\omega$ равны нулю: $\omega \times \omega=0,$ \, $\omega'\cdot\omega=0$.
Ранг матрицы $\mathrm{rot} (i\xi )$ равен двум при  $\xi \neq 0$,
поэтому  $\mathbf{w}=c\,\omega$, где $c$ - постоянная,  и других
решений нет.
Граничное условие приводит нас к уравнению:
$|\tau |c=0$.   Значит  $c=0$ при $|\tau| > 0$ и, следовательно, $\mathbf{v}=0$.

Итак, система $(2.4)$ с краевым условием $\mathbf{n}\cdot
\mathbf{u}{{|}_{\Gamma }}=g$  при $\lambda\neq 0$ - 
эллиптическая задача.
Мы скажем при этом, что задача (2.1)  при
$\lambda\neq 0$ является обобщенно эллиптической. 

Обобщенная эллиптичность задачи (2.2) доказана в  \cite{ds18}.

 \subsection{  Операторы задач (2.1) и (2.2) в пространствах ${\bf{H}^{s}}(G)$}
Пусть вектор-функция $\mathbf{u}$ принадлежит пространству Соболева
${\bf{H}^{s+1}}$, где $s\geq 0$--целое. Тогда компоненты
$\mathrm{rot} \mathbf{u}$ и $\mathrm{div} \mathbf{u}$ принадлежат
${H}^{s}(G)$, а вектор-функция $\mathbf{f}:=\mathrm{rot} \mathbf{u}+
\lambda \mathbf{u}$ принадлежит пространству
\begin{equation*}\klabel{pr 2} {\bf{E}^{s}}(G)=\{\mathbf{f}\in
{\mathbf{H}^{s}}: \mathrm{div} \mathbf{f}\in
{H}^{s}\},  \quad \|\mathbf{v}\|_{\mathbf{E}^{s}}=(
\|\mathbf{v}\|^2_{s }+
\|\mathrm{div}\mathbf{v}\|^2_{s})^{1/2}.  \eqno{(2.7)} \end{equation*}

Функция  $g:=\gamma({\mathbf{n}}\cdot\mathbf{u})\equiv\mathbf{n}
\cdot\mathbf{u}|_{\Gamma}$ принадлежит пространству
Соболева-Слободецкого $H^{s+1/2}(\Gamma)$. \newline Следовательно,
 задаче (2.1) соответствует ограниченный оператор
\begin{equation*}\klabel{op 1} \mathbb{A}\mathbf{u}\equiv\left( \begin{matrix}
\mathrm{rot} +\lambda I \\\gamma
\mathbf{n}\cdot \end{matrix}\right)\mathbf{u}:
\bf{H}^{s+1}(G)\rightarrow\left(
\begin{matrix}\mathbf{E}^{s}(G)\\ H^{s+1/2}(\Gamma)\end{matrix}\right).  \eqno{(2.8)}
\end{equation*}
Аналогично,  задаче (2.2) соответствует ограниченный оператор
\begin{equation*}\klabel{op 2} \mathbb{B}\mathbf{u}\equiv\left( \begin{matrix}
\nabla \mathrm{div} +\lambda I \\\gamma
\mathbf{n}\cdot \end{matrix}\right)\mathbf{u}:
\bf{H}^{s+2}(G)\rightarrow\left(
\begin{matrix}\bf{F}^{s}(G)\\ H^{s+3/2}(\Gamma)\end{matrix}\right),   \eqno{(2.9)}
\end{equation*}
 \[\bf{F}^{s}=\{\mathbf{f}\in
{\mathbf{H}^{s}}: \mathrm{rot} \mathbf{f}\in {H}^{s+1}\}, \quad
\|\mathbf{v}\|_{\mathbf{F}^{s}}=(
\|\mathbf{v}\|^2_{s }+
\|\mathrm{rot}\mathbf{v}\|^2_{s+1})^{1/2}.  \eqno{(2.10)}\]
Мы показали, что краевые задачи (2.1) и (2.2) являются обобщенно эллиптическими и что операторы  $\mathbb{A}$ и $ \mathbb{B}$
является {\it эллиптическими по Солонникову} [14].   \quad Из его Теоремы 1.1 следуют {\bf Теорема 1} и {\bf Теорема 2} (см. п.1.4).
 Область  $G$ ограничена с гладкой границей.
    \subsection{ Спектральные задачи  операторов  $ \text{rot}$ и  $\nabla\text{div}$}  
Они состоят в нахождении ненулевых вектор-функций ( полей)
$\mathbf{u}$ и $\mathbf{v}$ и чисел $\lambda$ и $\mu$  таких, что
\begin{equation*}
\klabel{srot__1_}
 \text{rot} \mathbf{u}=\lambda \mathbf{u}(\mathbf{x}), \quad
 \mathbf{x}\in G, \quad \gamma \mathbf{n}\cdot \mathbf{u}=0,
 \quad \mathbf{u}\in \mathcal{C}^1(G)\cap \mathcal{C}(\overline{G}),
 \eqno{(2.11)}
\end{equation*}
\begin{equation*}
\klabel{sgrd__1_}
 -\nabla\text{div}\mathbf{v}=\mu \mathbf{v}(\mathbf{x}), \quad
 \mathbf{x}\in G, \quad \gamma \mathbf{n}\cdot \mathbf{v}=0,
 \quad \mathbf{v}\in \mathcal{C}^2(G)\cap \mathcal{C}(\overline{G}).
\eqno{(2.12)}\end{equation*} Из теорем A, B и оценок вытекают полезные
свойства решений\newline {{\it спектральных задач} операторов ротора
и градиента дивергенции:\newline а)  {\it ненулевые собстенные
значения} имеет конечную кратность,
\newline б)  соответствующие  им {\it обобщенные собственные
функции} бесконечно дифференцируемы вплоть до границы области, то-есть
 поля $\mathbf{u}_{\lambda}(\mathbf{x})$ и
$\mathbf{v}_{\mu}(\mathbf{x}) \in  \mathcal{C}^{\infty}(\overline{G})$ при  $\lambda\ne 0$ и $\mu\ne 0$.

 {\it Замечание. } {\small   Мне удалось найти в \cite{saUMJ13} формулы  решений  спектральной задачи (2.11) в шаре   
   благодаря идее сведения задачи (2.1) к задаче Дирихле для уравнения Гельмгольца} {\small \footnote{Ее осуществил  мой выпусник НГУ 1971 года А.А.Фурсенко. В дипломной работе "Краевая задача для одной равномерно неэллиптической системы,"  \, 29 стр.,  он решил задачу (2.1)  в шаре в  классах Гельдера. } } и учебнику Владимирова \cite{vla}.
 	{\small   Поля
 	$\mathbf{u}^{\pm}_{\kappa}$, отвечающие
 	ненулевым  значениям  ротора	$\pm\lambda_{\kappa}=\pm\rho_{n,m}/R$,
   выражаются через 
 сферические функции и функции $\psi_n(z)$, см. (1.18) ,
 		где  $  \kappa =(n,m,k), \,  n,  \,m\in \mathbb{N},\,   |k|\leq n$,  а
 	числа $\pm\rho_{n,m}$ - нули функций  $\psi_n (r)$. 	
 	Поля  $\mathbf{q}_{\kappa}$  	со значениями $-\nu_{\kappa}^2$,  где
 	$\nu_{\kappa}=\alpha_{n,m}/R$, определяются решениями   задачи Неймана;   $\alpha_{n,m}$ - нули    производных 
 	$\psi'_n(r)$, $ n\geq 0 $ .}  
 	
 		Поля $\{\mathbf{u}^{+}_{\kappa}\}\cup\{\mathbf{u}^{-}_{\kappa}\}\cup \{\mathbf{q}_{\kappa}\}$	образуют базис в  	 $\mathbf{L}_{2}(B)$.} 
\section{Класс  $\mathbf{V}^{0}$ и его подпространства $\mathbf{W}^{k}$}

Другой путь решения задачи (2.1) открылся после работы Joshida и Giga \cite{yogi}. 
Они ввели оператор $S:\mathbf{V}^{0}\rightarrow   \mathbf{V}^{0}$  в пространстве  $L^2_{\Sigma}\equiv\mathbf{V}^{0}$  с областью определения ${H}^{1}_{\Sigma\Sigma}\equiv \mathbf{W}^{1}$,  который совпадает с $\mathrm{rot}\mathbf{u}$, если   $\mathbf{u}\in \mathbf{W}^{1}$,
и доказали теорему (цитирую):
 {\it The operator $S$ is self-adjoint in the space $L^2_{\Sigma}$.
 The spectrum $\sigma(S)$ of $S$ consists of only point
 spectrum $\sigma_p(S)\subset\mathbb{R}$. Therefore, the set of
 eigenfunctions of $S$ gives an orthogonal complete basis of
 the space $L^2_{\Sigma}$.}

 \subsection { Собственные поля ротора и ортогональный 	базис в  	 $ \mathbf{V}^{0}$ }  Поля $\mathbf{u}_{\lambda}(\mathbf{x})$ 
  принадлежат пространствам $\mathbf{W}^1(G)\cap \mathcal{C}^\infty(\overline{G})$.  Из соотношения 
$ (\mathrm{rot}+\lambda {I})
 (\mathrm{rot}- \lambda { I})\mathbf{ u}=-
 \Delta \mathbf{u}
 +\nabla \mathrm{div} \mathbf{ u}-\lambda^2 \mathbf{ u}$ 
и определения пространства $V^0(G)$
 видим, что собственные поля ротора
 $\mathbf{u}^{\pm}_{\lambda}$,
 отвечающие  значениям $\pm\lambda\ne 0 $,
 является также
 собственными полями оператора Лапласа: 
 \begin{equation*}\klabel{ladi__2_}
  -\Delta \mathbf{u}=\lambda^{2}\mathbf{u},\quad
\mathrm{div} \mathbf{u}=0, \,\mathbf{n}\cdot \mathbf{u}|_{\Gamma}=0.
\eqno{(3.1)}\end{equation*}

 {\small Множество собственных значений $\mu=\lambda^2$ этого оператора 
 счетно, положительно и каждое из них имеет
 конечную кратность. Перенумеруем их в порядке возрастания:
 $0<\mu_1\leq \mu_2\leq ...$, повторяя $ \mu_k$ столько раз, какова
 его кратность. Соответствующие вектор-  функции  обозначим через
$\mathbf{u}_{1}^{\pm}, \mathbf{u}_{2}^{\pm}$, ..., так чтобы каждому
   значению $\pm\sqrt{\mu_{k}}$ соответствовала только
одна  функция $\mathbf{u}_{k}^{\pm}$: $\mathrm{rot}
\mathbf{u}_{k}^{\pm}=\pm\sqrt{\mu_k} \mathbf{u}_{k}^{\pm}$,\quad
$k=1,2,...$.

 Собственные функции, соответствующие одному и тому же
собственному значению, выберем ортонормальными, используя процесс
ортогонализации Шмидта  (см. \cite{vla}). 
 Поля, соответствующие
различным с.- значениям, ортогональны. Их нормируем.
Нормированные собственные поля ротора обозначим через
$\mathbf{q}^{\pm}_{j}$,  норма  $\|\mathbf{q}^{\pm}_{j}\|=1$.
 Они составляют полный ортонормированный базис в классе 
$\mathbf{V}^0$ вихревых полей в $\mathbf{L}_2(G)$.}

 \subsection{  Ряды Фурье в  $\mathbf{V}^0$ }
Проекция вектор-функции $\mathbf{f}$ из $\mathbf{L}_2(G)$ на
$\mathbf{V}^0$ имеет вид:
\begin{equation*}
\klabel{sp__1_}\mathbf{f}_{\mathbf{v}}=\sum_{j=1}^{\infty}[(\mathbf{f},\mathbf{q}^{+}_{j})\mathbf{q}^{+}_{j}+
(\mathbf{f},\mathbf{q}^{-}_{j})\mathbf{q}^{-}_{j}].\eqno{(3.2)}
\end{equation*}
Действительно, частичные суммы $\mathbf{f}^n_{\mathbf{v}}$ этого
ряда состоят из элементов, для которых $0<\lambda_j\leq N(n)$:
\[\mathbf{f}^n_{\mathbf{v}}=\sum_{j=1}^{n}[(\mathbf{f},\mathbf{q}^{+}_{j})\mathbf{q}^{+}_{j}+
(\mathbf{f},\mathbf{q}^{-}_{j})\mathbf{q}^{-}_{j}],  \quad \|\mathbf{f}^n_{\mathbf{V}}\|^2
\leq
\|\mathbf{f}\|^2 ,    \eqno{(3.3)}\] проекции
$(\mathbf{f}-\mathbf{f}^n_{\mathbf{v}},\mathbf{q}^{\pm}_{j} )=0$,
если 
 \, \quad
$0<\lambda_j\leq  N(n) $,\quad и
\[\|\mathbf{f}_{\mathbf{v}}-\mathbf{f}^n_{\mathbf{v}}\|^2=
\|\mathbf{f}_{\mathbf{v}}\|^2-\|\mathbf{f}^n_{\mathbf{v}}\|^2
\rightarrow 0 \quad \text{при} \quad n\rightarrow\infty.\] По построению  $\mathbf{f}^n_{\mathbf{v}}\in
\mathcal{C}^{\infty}(\overline{G})$, 
$\mathrm{div}\mathbf{f}^n_{\mathbf{v}}=0$,
$\gamma_\mathbf{n}\mathbf{f}^n_\mathbf{V}=0$ и при любом $n$ поле
$\mathbf{f}^n_{\mathbf{V}}\bot Ker\mathbf{(rot)}$ в
$\mathbf{L}_2(G)$.
 Значит,
 $(\mathbf{f}^n_{\mathbf{V}}, \nabla h)=0$ для любой функции
 $h\in H^1(G)$. Переходя к пределу
 получим
 $(\mathbf{f}_{\mathbf{V}}, \nabla h)=0$, то-есть
 вектор $\mathbf{f}_{\mathbf{V}}\perp \mathcal{A}\subset Ker(\mathbf{rot})$. 
 Он принадлежит
 $\mathbf{V}^0$, если пространство
 $\mathcal{B}_H\subset Ker(\mathbf{rot})$
 пусто.
 В общем случае,
 \begin{equation*}\klabel{obh__1_}
 \mathbf{f}_{\mathbf{V}}\in  \mathbf{V}^0 \quad  \Leftrightarrow
 \quad 
 (\mathbf{f}_{\mathbf{V}},\mathbf{h}_{i})=
 0\quad \forall\,   i=1,...,\rho.   .\eqno{(3.4)} \end{equation*}
 Так как 
$ \mathrm{rot} (\mathbf{f}^n_{\mathbf{V}})=
 \sum_{j=1}^{n}\lambda_j
[(\mathbf{f},\mathbf{q}^{+}_{j})\mathbf{q}^{+}_{j}-
(\mathbf{f},\mathbf{q}^{-}_{j})\mathbf{q}^{-}_{j}]$
  и суммы
  $\mathbf{f}^n_{\mathbf{V}}$\,\,и $\mathrm{rot}
(\mathbf{f}^n_{\mathbf{V}})$ принадлежат $\mathbf{V}^{0}$, то
$\mathbf{f}^n_{\mathbf{V}}\in\mathbf{W}^{1}$
- области определения оператора $S$.

 {\it По
определению $S\mathbf{w}=\mathrm{rot}\mathbf{w}$ для любого
$\mathbf{w}\in \mathbf{W}^{1}$.} Следовательно
\begin{equation*}\klabel{ops__1_} S\mathbf{f}_{\mathbf{V}}=\lim_{n\rightarrow\infty}
\mathrm{rot}
(\mathbf{f}^n_{\mathbf{V}})=\sum_{j=1}^{\infty}\lambda_j
[(\mathbf{f},\mathbf{q}^{+}_{j})\mathbf{q}^{+}_{j}-
(\mathbf{f},\mathbf{q}^{-}_{j})\mathbf{q}^{-}_{j}],  \eqno{(3.5)}\end{equation*}
 если ряд
сходится и принадлежит $\mathbf{V}^{0}$.  Ясно,  что     $S\mathbf{f}_{\mathbf{V}}\in\mathbf{V}^{0}$,   если
 \begin{equation*}\klabel{obh__2_}f\in\mathbf{H}^1(G),
 \quad (\mathbf{f}_{\mathbf{V}},\mathbf{h}_{i})=0\quad \text{и}\quad
(S\mathbf{f}_{\mathbf{V}},\mathbf{h}_{i})=0 \quad
 \forall \, i=1,...,\rho.   \eqno{(3.6)}\end{equation*}
Легко видеть, что оператор $S$ замкнут. 
Следовательно, предел \newline 
  $S\mathbf{f}_{\mathbf{V}}$ {\it не зависит
от  выбора в $\mathbf{V}^0$ последовательности  $ \mathbf{w}_n \rightarrow \mathbf{f}_{\mathbf{V}}$.}
 \subsection{Подпространства $\mathbf{V}^0$} Введем пространства
 \begin{equation*}\klabel{prow__k_} \mathbf{W}^{k}(G)= 
\{\mathbf{f}\in\mathbf{V}^{0}(G),..., \mathrm{rot}^k
\mathbf{f}\in\mathbf{V}^{0}(G) \}\quad \forall\, k\geq
1.     \eqno{(3.7)}\end{equation*}
Вложение
$\mathbf{W}^1\subset\mathbf{H}^{1}(G)$  вытекает из оценки  (1.4) при $s=0$:
\[C_0\|\mathbf{f} \|_1\leq \| \mathrm{rot} \mathbf{f} \|+\|\mathbf{f} \| , \quad  C_0>0.     \eqno{(3.8)}\]
 По индукции
$\mathbf{W}^k\subset\mathbf{H}^{k}(G)$. \quad 
 Очевидно, что $\mathbf{W}^k\subset...\subset \mathbf{W}^1\subset \mathbf{V}^{0}.$
\newline При  $n<\infty$ ряды  $\mathbf{f}^n_{\mathbf{V}}$  принадлежат любому из этих пространств. \newline 
  Оператор
$S$ отображает  $\mathbf{W}^{k}$ на $\mathbf{W}^{k-1}$  при  $k>1$,  а  $\mathbf{W}^{1}$ на $\mathbf{W}^{0}\equiv
  \mathbf{V}^{0}$.

 Пространство $\mathbf{V}^{0}$ ортогонально ядру ротора
 в $\mathbf{L}_{2}(G)$, поэтому $S$ имеет единственный обратный
 оператор $S^{-1}$, определенный на $\mathbf{V}^{0}$:
\begin{equation*}\klabel{obr__2_}S^{-1}\mathbf{f}_{\mathbf{V}}=\lim_{n\rightarrow\infty}
S^{-1} (\mathbf{f}^{n}_{\mathbf{V}})=
\sum_{j=1}^{\infty}\lambda_j^{-1}
[(\mathbf{f},\mathbf{q}^{+}_{j})\mathbf{q}^{+}_{j}-
(\mathbf{f},\mathbf{q}^{-}_{j})\mathbf{q}^{-}_{j}]. \eqno{(3.9)} \end{equation*}

Доказано   \cite{ yogi}, что оператор $S^{-1}$ - компактен.

 Следствие. {\it Спектр оператора $S^{-1}$ точечный с единственной точкой
 накопления в нуле,\quad
 $\lambda^{-1}_j\rightarrow 0$ при ${j\rightarrow\infty}$.}

 Очевидно, что оператор $S^{-1}:\mathbf{V}^{0}\rightarrow \mathbf{W}^{1}$ 
 и 
  $S^{-1}:\mathbf{W}^{k-1}\rightarrow \mathbf{W}^{k}$.

\subsection{Полнота пространства $\mathbf{V}^{0}$} В базисе из собственных
функций ротора скалярное произведение векторов
$\mathbf{f},\mathbf{g}\in \mathbf{V}^{0}$ имеет вид:
\begin{equation*} \klabel{spb 1}
( \mathbf{f},\mathbf{g})= \lim_{n\rightarrow\infty} (
\mathbf{f}_{\mathbf{v}}^n, \mathbf{g}_{\mathbf{v}}^n)=
\sum_{j=1}^{\infty}
[(\mathbf{f},\mathbf{q}_{j}^+)(\mathbf{g},\mathbf{q}_{j}^+)
+(\mathbf{f},\mathbf{q}_{j}^-)(\mathbf{g},\mathbf{q}_{j}^-)]. \eqno{(3.10)}
\end{equation*}
Согласно  Владимирову \cite{vla},  ортонормальная система $\{\mathbf{q}_{j}^+\}\cup \{
\mathbf{q}_{j}^-\}_{j=1,2,...}$ полна в $\mathbf{V}^{0}$, если для
любой $\mathbf{f}\in \mathbf{V}^{0}$ ее ряд  (3.2) сходится к
$\mathbf{f}$ в $\mathbf{L}_{2}(G)$. По Теореме 1 $\S 1.9$ эта система полна в
$\mathbf{V}^{0}$, тогда и только тогда, когда для любой функции
$\mathbf{f}\in \mathbf{V}^{0}$ выполняется равенство
Парсеваля-Стеклова, которое называется уравнением замкнутости:
\begin{equation*} \klabel{UrZa 1 }
\sum_{j=1}^{\infty}[(\mathbf{f},
\mathbf{q}^{+}_{j})^2+(\mathbf{f},\mathbf{q}^{-}_{j})^2]=
\|\mathbf{f}\|^2. \eqno{(3.11)}
\end{equation*}
 Пространство $\mathbf{W}^1$ плотно
 в $\mathbf{V}^{0}$, так как множество
 $\mathbf{C}_0^{\infty}(G) \cap\mathbf{V}^{0}$,  плотное в  $\mathbf{V}^{0}$, содержится
 в $\mathbf{W}^1$.
 Квадрат нормы   $\mathbf{f}\in\mathbf{C}_0^{\infty}(G)\cap\mathbf{V}^{0}$ ограничен:
\[\|\mathbf{f}\|^2_{\mathbf{W}^1}=\|\mathbf{f}\|^2+
\|\mathrm{rot}\mathbf{f}\|^2=\sum_{j=1}^{\infty}
(1+{\lambda}_{j}^2)[(\mathbf{f},
\mathbf{q}^{+}_{j})^2+(\mathbf{f},\mathbf{q}^{-}_{j})^2]<\infty
\Rightarrow \]
\[\lim_{n\rightarrow\infty} \|\mathbf{f}^n_{\mathbf{V}}\|^2=
\sum_{j=1}^{\infty}[(\mathbf{f},
\mathbf{q}^{+}_{j})^2+(\mathbf{f},\mathbf{q}^{-}_{j})^2]=
\|\mathbf{f}\|^2,\quad  \text{ч.т.д.} .\]

\subsection{ Оператор $S:\mathbf{V}^0\rightarrow \mathbf{V}^0$ - самосопряжен }  Действительно,
если $\mathbf{f}$ и $\mathbf{g}$ принадлежат $\mathbf{W}^1$, то
имеют место равенства
\begin{equation*} \klabel{sor }(\text{S} \mathbf{f}, \mathbf{g})=
(\mathbf{f},\text{S}\mathbf{g})= \sum_{j=1}^{\infty}
{\lambda}_{j}[(\mathbf{f},\mathbf{q}_{j}^+)(\mathbf{g},\mathbf{q}_{j}^+)
-(\mathbf{f},\mathbf{q}_{j}^-)(\mathbf{g},\mathbf{q}_{j}^-) ]. \eqno{(3.12)}
\end{equation*}
 Отметим, что что равенство \quad 
$\klabel{sic__1_} \int_G (\mathrm{rot} \mathbf{u})\cdot \mathbf{v}
d\mathbf{x}=\int_{G}\mathbf{u}\cdot (\mathrm{rot} \mathbf{v}) d
\mathbf{x}$ 
для любых функций $\mathbf{u}$ и $\mathbf{v}$
из $\mathcal{D}(S)$ доказано в \cite{yogi}, а в случае шара - в  \cite{saUMJ13},
другим способом.
 Доказано также  \cite{yogi},  что оператор  $S$ самосопряжен.

 \subsection{Оператор $S+\lambda I:\mathbf{W}^{1}
	\rightarrow \mathbf{V}^0$ - фредгольмов}
Действительно, по определению оператор $ S+\lambda I$
совпадает с
$\mathrm{rot}+\lambda I$   на $\mathbf{W}^{1}$.
При $\mathbf{f}\in \mathbf{W}^{1}$
\[(S+\lambda I)\mathbf{f}=\lim_{n\rightarrow\infty}
(\mathrm{rot+\lambda I}) (\mathbf{f}^n_{\mathbf{V}})
=\sum_{j=1}^\infty[(\lambda+\lambda_j)(\mathbf{f},\mathbf{q}^{+}_{j})\mathbf{q}^{+}_{j}+
(\lambda-\lambda_j)(\mathbf{f},\mathbf{q}^{-}_{j})\mathbf{q}^{-}_{j}]\]
и ряд сходится в $\mathbf{L}_{2}(G)$,  так как \[\|  (S+\lambda I)\mathbf{f}\|^2_{ \mathbf{V}^0}=\sum_{j=1}^\infty[|\lambda+\lambda_j|^2(\mathbf{f},\mathbf{q}^{+}_{j})^2+|\lambda-\lambda_j|^2(\mathbf{f},\mathbf{q}^{-}_{j})^2]\leq \]
\[c^2_0\sum_{j=1}^{\infty}
(1+{\lambda}_{j}^2)[(\mathbf{f},
\mathbf{q}^{+}_{j})^2+(\mathbf{f},\mathbf{q}^{-}_{j})^2]= c^2_0 \|\mathbf{f}\|^2_{\mathbf{W}^1}, \quad   \eqno{(3.13)}\]
 \[c^2_0 = max _j \, (a_j^+, a_j^-)\,   \text{и} \quad a_j^{\pm }= (|1\pm \lambda/\lambda_j|^2/(1+1/{\lambda}_{j}^2))<\infty,  \,\,  \eqno{(3.14)} \]
 так  как  при больших $\lambda_{j}$ они находятся в окрестности 1.

Обратный оператор имеет вид:
\begin{equation*}\klabel{sp__3_}
(S+\lambda I)^{-1}\mathbf{f}=
\sum_{j=1}^{\infty}\left[
\frac{(\mathbf{f},{\mathbf{q}}_{j}^+)}{\lambda+\lambda_{j}}
\mathbf{q}_{j}^+(\mathbf{x})+
\frac{(\mathbf{f},{\mathbf{q}}_{j}^-)}{\lambda-\lambda_{j}}
\mathbf{q}_{j}^-(\mathbf{x})\right], \eqno{(3.15)}\end{equation*}
если ни одно из слагаемых этого ряда не обращается в
бесконечность. Это означает, что либо $\lambda\pm \lambda_{j}\neq 0$ для всех $j$, либо  $(\mathbf{f},\mathbf{q}^{-}_{j})=0$, если  $\lambda=\lambda_j=\lambda_{j_0}$,  и  функция $\mathbf{f}$
ортогональна всем собственным полям
$\mathbf{q}^{-}_{j}(\mathbf{x})$ ротора, отвечающим собственному
значению $\lambda_{j_0}$. При этом 
\[  \|(S+\lambda I)^{-1}\mathbf{f}\|^2_{\mathbf{W}^1}=\sum_{j=1}^{\infty}\left[
\frac{(1+ \lambda_{j}^2)}{|\lambda+\lambda_{j}|^2}(\mathbf{f},{\mathbf{q}}_{j}^+)^2+\frac{(1+ \lambda_{j}^2)}{|\lambda-\lambda_{j}|^2}(\mathbf{f},{\mathbf{q}}_{j}^-)^2\right] \leq
C^2_0 \|\mathbf{f}\|^2_{\mathbf{V}^0}, \] 
\[C^2_0 = max _j \, (A_j^+, A_j^-) \,\,\text{ и } A_j^{\pm }= ((1+1/{\lambda}_{j}^2) /|1\pm \lambda/\lambda_j|^2<\infty.    \eqno{(3.16)}\] 
Итак, оба оператора непрерывны и имеет место {\bf Теорема 3}  (см. п.1. 5).

 Аналогично предыдущему видим, что \quad  $ \|  (S+\lambda I)\mathbf{f}\|^2_{ \mathbf{W}^m}=$
\[\sum_{j=1}^\infty (1+ \lambda_{j}^{2m})[|\lambda+\lambda_j|^2(\mathbf{f},\mathbf{q}^{+}_{j})^2+|\lambda-\lambda_j|^2(\mathbf{f},\mathbf{q}^{-}_{j})^2]\leq  \]
\[c^2_m\sum_{j=1}^{\infty}
(1+{\lambda}_{j}^{2(m+1)})[(\mathbf{f},
\mathbf{q}^{+}_{j})^2+(\mathbf{f},\mathbf{q}^{-}_{j})^2]= c^2_m \|\mathbf{f}\|^2_{\mathbf{W}^{m+1}},    \quad   \eqno{(3.17)}\]
\[  \|(S+\lambda I)^{-1}\mathbf{f}\|^2_{\mathbf{W}^{m+1}}\leq  C^2_m \|\mathbf{f}\|^2_{\mathbf{W}^{m}}, \quad   c_m,\,  C_m< \infty .  \quad   \eqno{(3.18)} \]

Числа $ c_m$ и  $C_m$ при $m=0$ совпадают с  (3.14)  и  (3.16).  
 По определению  $\mathbf{W}^{0}\equiv \mathbf{V}^{0}$.\quad  
 Из  Теоремы 3 и оценок следует {\bf Лемма 1} п.1.5 о свойствах отображений  $(S+\lambda I)$ и  $(S+\lambda I)^{-1}$. 
Так же доказана  {\bf Лемма 2} п.1.7 о свойствах отображений  $( \mathcal{N}_d+\lambda I)$ и  $( \mathcal{N}_d+\lambda I)^{-1}$.

\subsection{Сходимость ряда Фурье в норме пространства $\mathbf{H}^k(\Omega)$}
 Область
$\Omega$, гомеоморфную шару, выделил С.Л.Соболев в \cite{sob54}. В этом случае  пространство $\mathcal{B}_H(\Omega)$  пусто. Границу  области $\Omega$  будем предполагать гладкой.
Скалярное произведение
в $\mathbf{H}^k(\Omega)$ Сергей  Львович определяет так:
\begin{equation*} \klabel{skpro 1}
(\mathbf{f},\mathbf{g})_k=(\mathbf{f},\mathbf{g})+
 \int_{\Omega} \sum_{|\alpha|=k}\frac{k!}{\alpha
 !}\partial^{\alpha}\mathbf{f}\cdot\partial^{\alpha}\mathbf{g}
 d \mathbf{x},\quad k\geq 1.  \eqno{(3.19)}
\end{equation*}
 В  пространстве $ \mathbf{W}^k(\Omega)=
\{\mathbf{f}\in \mathbf{V}^0,
..., \mathrm{rot}^{k} \mathbf{f}\in \mathbf{V}^0\}$
норму $\mathbf{f}\in{\mathbf{W}^k}$ выберем   так же:
$\|\mathbf{f}\|^2_{\mathbf{W}^k}\equiv \|\mathbf{f}\|^2+\|\mathrm{rot}^k \mathbf{f}\|^2$.  \quad

 Имеет место Теорема 4 (см. п.1.6).

Доказательство Теоремы 4. Граница $\partial\Omega\in \mathcal{C}^\infty$ и
собственные функции $\mathbf{q}_{j}^{\pm}(\mathbf{x})$ оператора ротор
принадлежат классу
$\mathcal{C}^{\infty}(\overline{\Omega})$,
а, значит, любому из пространств $\mathbf{W}^l(\Omega)$, $l>0$.
Поэтому, если ряд   Фурье (1.11) вектор-функции
$\mathbf{f}$  из
$\mathbf{H}^k(\Omega)$ сходится в норме
$\mathbf{H}^k(\Omega)$, то $\mathbf{f}\in\mathbf{V}^0,...,
\mathrm{rot}^{k} \mathbf{f}\in
\mathbf{V}^0\subset\mathbf{L}_2(\Omega)$ и, значит,
$\mathbf{f}\in\mathbf{W}^k(\Omega)$.\quad
Необходимость доказана.

Пусть $\mathbf{f}\in \mathbf{W}^k(\Omega)$, где $k\geq 1$.
 Приведем доказательство   неравенства (1.12).
Так как $\mathcal{S}{\mathbf{f}}=\text{rot} {\mathbf{f}}$ на
$\mathbf{W}^k\subset \mathbf{W}^1(\Omega)$
и $\mathcal{S} \mathbf{q}^{\pm}_{j}={\pm}\lambda_j
\mathbf{q}^{\pm}_{j}$, то
\begin{equation*}
\klabel{rot 3} (\text{rot} {\mathbf{f}},\mathbf{q}^{\pm}_{j}) =
{\pm}\lambda_j({\mathbf{f}} ,\mathbf{q}^{\pm}_{j}).\eqno{(3.20)}
\end{equation*}
Пусть $\beta_{k,j}^{\pm}$
коэффициенты Фурье функции $\text{rot}^k \mathbf{f}$. По формуле
(3.20)
\begin{equation*}
\klabel{rot k}\beta_{k,j}^{\pm}=(\text{rot}^k \mathbf{f},
\mathbf{q}^{\pm}_{j})={\pm}\lambda_j(\text{rot}^{k-1} {\mathbf{f}}
,\mathbf{q}^{\pm}_{j}) =...= ({{\pm}\lambda_j})^{k}( {\mathbf{f}}
,\mathbf{q}^{\pm}_{j}).     \eqno{(3.21)} \end{equation*}

Поскольку $\text{rot}^k \mathbf{f}\in
\mathbf{L}_2(\Omega)$, то  $\sum_{j=1}^{\infty} [(\beta_{k,j}^+)^2+
(\beta_{k,j}^-)^2]= \|\text{rot}^k \mathbf{f}\|^2.$

Итак, для вектор-функций $\mathbf{f}\in \mathbf{W}^k(\Omega)$ имеем
\begin{equation*} \klabel{orf 4}
\sum_{j=1}^{\infty} {\lambda}_{j}^{2k}
((\mathbf{f},\mathbf{q}_{j}^+)^2 +(\mathbf{f},\mathbf{q}_{j}^-)^2)=
\|\text{rot}^k \mathbf{f}\|^2\leq C
\|\mathbf{f}\|^2_{\mathbf{H}^k(\Omega)}. \eqno{(3.22)}
\end{equation*}
Последнее неравенство вытекает из определений нормы в $\mathbf{H}^k(\Omega)$.
Неравенство (1.12) доказано. Докажем сходимость ряда
(1.11) к $\mathbf{f}$ в  норме $\mathbf{H}^k(B)$.
Пусть $\mathbf{S}_l(\mathbf{x})$ - частичная сумма ряда
(1.11). Очевидно, что $\mathbf{S}_l(\mathbf{x})\in
\mathbf{W}^l(\Omega)$ при $l>0$. В частности, $\text{div}
\mathbf{S}_l(\mathbf{x})=0$ и
$\gamma\mathbf{n}\cdot\mathbf{S}_l(\mathbf{x})=0$.
Поэтому оценка
(1.4) при $s=0$ принимает вид \quad $ C_1\|\mathbf{S}_l\|_{1}
\leq\|\mathrm{rot} \mathbf{S}_l\|+ \|\mathbf{S}_l\|.$
Поскольку $\lambda_{j}^{-2}\rightarrow 0$ при $j\rightarrow \infty$,
норма $\|\mathbf{S}_l\|^2\leq c\|\mathrm{rot}\mathbf{S}_l\|^2$, где
$c=max_{j}\lambda_{j}^{-2}$.
Поэтому
$ \|\mathbf{S}_l\|^2_{1} \leq a_1\|\mathrm{rot} \mathbf{S}_l\|^2$ и
по индукции $ \|\mathbf{S}_l\|^2_{k}\leq a_k\|\mathrm{rot}^k \mathbf{S}_l\|^2$.

Пусть $\mathbf{f}\in \mathbf{W}^k(\Omega)$, где $k>0$. Согласно
неравенству (1.12), ряды в его левой части
сходятся и если $l>m\geq 1$, то
\[\|\mathbf{S}_l-\mathbf{S}_m\|^2_{k}\leq a_k\|rot^k
(\mathbf{S}_l-\mathbf{S}_m)\|^2
\leq \,a_k\sum_{m+1}^l\lambda_{j}^{2k}(|(\mathbf{f},\mathbf{q}_{j}^+)|^2
+|(\mathbf{f},\mathbf{q}_{j}^-|^2))\rightarrow  0\]  при
$l,m\rightarrow\infty$. Это означает, что ряд (1.11) сходится
к $\mathbf{f}$ в  норме $\mathbf{H}^k(B)$.

{\bf Замечание.}  Известно
вложение пространств
$\mathbf{H}^k(\Omega)\subset\mathbf{C}^{k-2}(\overline{\Omega})$ при
$k\geq 2$ в трехмерной области $\Omega$ и оценка:\quad  $\|\mathbf{f}\|_{\mathbf{C}^{k-2}(\overline{\Omega})}\leq
C_k\|\mathbf{f}\|_{\mathbf{H}^k (\Omega)}$\quad  для любой функции
$\mathbf{f}\in \mathbf{H}^k (\Omega)$, причем, постоянная $C_k>0$ не
зависит от $\mathbf{f}$ (см., например, Теорему 3 $\S 6.2$ в \cite{mi}).
В частности,
\begin{equation*}
\klabel{ck l}\|\mathbf{S}_l-\mathbf{S}_m\|_{\mathbf{C}^{k-2}(\overline{\Omega})}
\leq C_k\|\mathbf{S}_l-\mathbf{S}_m\|_{\mathbf{H}^k (\Omega)}.  \eqno{(3.23)}\end{equation*} 
Если
$\|\mathbf{S}_l-\mathbf{S}_m\|_{\mathbf{H}^k (\Omega)}\rightarrow 0$ при
$l,m\rightarrow\infty$, то
$\|\mathbf{S}_l-\mathbf{S}_m\|_{\mathbf{C}^{k-2}(\overline{\Omega})}\rightarrow
0$. Это означает, что ряд (1.11) сходится к $\mathbf{f}$ в
$\mathbf{C}^{k-2}(\overline{\Omega})$. Теорема доказана.

\section{Краевые задачи  в $\mathbf{L}_{2}(\Omega)$}

\subsection {Классы $C(2k, m)$ подпространств в  $\mathbf{L}_2(\Omega)$ }  Предположим, что область $\Omega$  гомеоморфна шару.  Если собственные поля $ \mathbf{q}_{j}(\mathbf{x})$  и  $\mathbf{q}_{j}^\pm(\mathbf{x})$  градиента дивергенции и ротора известны, то  элементы    $\mathbf{f}_\mathbf{A}\in \mathcal{A}$   и  $\mathbf{f}_\mathbf{V}\in \mathcal{B}=\mathbf{V}^0 $   
представляются   рядами Фурье:
\[\mathbf{f}_\mathcal{A}=
\sum_{j=1}^{\infty} (\mathbf{f},{\mathbf{q}}_{j}) \mathbf{q}_{j }(\mathbf{x}), \,\, 	\mathbf{f}_\mathbf{V}=\sum_{j=1}^{\infty}\left[(\mathbf{f},{\mathbf{q}}_{j}^{+})
 \mathbf{q}_{j}^{+}(\mathbf{x})+
(\mathbf{f},{\mathbf{q}}_{j}^{-})
\mathbf{q}_{j}^{-}(\mathbf{x})\right], \eqno{(4.1)}\]
а   элемент    $\mathbf{f}$  из $\mathbf{L}_2(\Omega)$ -- их суммой $\mathbf{f}_\mathcal{A}+\mathbf{f}_\mathbf{V}$.  Причем $ \text{div}\,\mathbf{f}=  \text{div}\,\mathbf{f}_\mathbf{A}$,  а   $\text{rot}\,\mathbf{f}=  \text{rot}\,\mathbf{f}_\mathbf{V}$,  так как $  \text{rot}\,\mathbf{f}_\mathbf{A}=0$  в $\mathcal{A}$ и $\text{div}\,\mathbf{f}_\mathbf{V}=0$   в $\mathcal{B}$.
Скалярное произведение  $(\mathbf{f}, \mathbf{g})$ полей
 $\mathbf{f}, \mathbf{g}$  из $\mathbf{L}_2(\Omega)$ равно   $(\mathbf{f}_\mathcal{A}, \mathbf{g}_\mathcal{A})+(\mathbf{f}_\mathbf{V}, \mathbf{g}_\mathbf{V})$.
Представления
 операторов    $\mathcal{N}_d$  в $ \mathcal{A}_{\gamma }$,  $S$ в  	
 	$ \mathcal{B}$ и  обратных  имеют вид
\[  \mathcal{N}_d \mathbf{f}_\mathcal{A}=
-\sum_{j=1}^{\infty}\nu ^2_j (\mathbf{f},{\mathbf{q}}_{j}) \mathbf{q}_{j }, \quad  	S \mathbf{f}_\mathbf{V}=\sum_{j=1}^{\infty}\lambda_{j}\left[(\mathbf{f},{\mathbf{q}}_{j}^{+})
\mathbf{q}_{j}^{+}-
(\mathbf{f},{\mathbf{q}}_{j}^{-})
\mathbf{q}_{j}^{-}\right], \]
 \[  \mathcal{N}_d  ^{-1}\mathbf{f}_\mathcal{A}=
 -\sum_{j=1}^{\infty}\nu ^{-2}_j (\mathbf{f},{\mathbf{q}}_{j}) \mathbf{q}_{j },\, \,  	S^{-1}\mathbf{f}_\mathbf{V}=\sum_{j=1}^{\infty}\lambda_{j}^{-1}\left[(\mathbf{f},{\mathbf{q}}_{j}^{+}) \mathbf{q}_{j}^{+}-
 (\mathbf{f},{\mathbf{q}}_{j}^{-}) \mathbf{q}_{j}^{-} \right]. \]
 
Рассмотрим пространства
\[ \mathbf{A}^{2k}_{\gamma }\equiv \{\mathbf{f}\in \mathcal{A}_{\gamma },...,(\nabla  \mathrm{div})^k \mathbf{f}\in \mathcal{A}_{\gamma }\}\,\text{и}\,  \mathbf{W}^m \equiv \{\mathbf{g}\in \mathbf{V}^0,..., (\mathrm{rot})^m \mathbf{g} \in \mathbf{V}^0\} \]  при  $k\geq 1, m\geq1$;
$ \mathbf{A}^{0}_{\gamma }\equiv  \mathcal{A},   \mathbf{W}^{0}\equiv  \mathbf{V}^{0}\equiv   \mathcal{B}$. 

 Имеют место вложения
 \[  \mathbf{A}^{2k}_{\gamma }\subset  \mathbf{A}^{2(k-1)}_{\gamma }  \subset  \mathbf{A}^{2}_{\gamma } \subset  \mathcal{A}_{\gamma } , \quad  \mathbf{W}^m\subset \mathbf{W}^{m-1}\subset...\subset \mathbf{W}^1\subset \mathbf{V}^{0}.\eqno{(4.2)}\]
Прямую сумму векторных пространств  $ \mathbf{A}^{2k}_{\gamma }\oplus  \mathbf{W}^m$  обозначим  как $C(2k, m)$    и назовем классом; $k\geq 0, m\geq 0 $  -целые, $k+m>0$.  
   Операторы $(\mathcal{N}_d^{-1}, I)$, $(I, S^{-1})$  и $(\mathcal{N}_d^{-1}, S^{-1})$ отображают  класс  $C(2k, m)$
на классы $C(2(k+1), m)$,  $C(2k, m+1)$ и  $C(2(k+1), m+1)$ и обратно
 (пп. 1.5 и 1.7).

  \subsection{ Пространство   $ \mathbf{E}^{s}(\Omega)$, $s\geq 0$ целое} Оно определяется  в \cite{rt} так     
  $ \mathbf{E}^{s}=\{\mathbf{f}\in
\mathbf{H}^{s}: \text{div}\mathbf{f}\in
{H}^{s}\}$. Квадат  нормы $  \| \mathbf{f}\|_{E^s}^2= \| \mathbf{f}\|_s^2+ \| \text{div}\mathbf{f}\|_s^2 = \| \mathbf{f}_\mathcal{A}\|_s^2+ \|\text{div} \mathbf{f}_\mathcal{A}\|_s^2+ \| \mathbf{f}_\mathbf{V}\|_s^2 .$\quad $ \mathbf{E}^{s}$ -
 пространство Гильберта и
\begin{equation*} \klabel{vlo s}
\mathbf{C}_0^{\infty}(\Omega)\subset \mathbf{E}^{s}(\Omega), \quad
\mathbf{H}^{s+1}(\Omega)\subset \mathbf{E}^{s}(\Omega)
\subset\mathbf{H}^s(\Omega). \eqno{(4.3)}\end{equation*}

Очевидно, что  ${\bf rot} \mathbf{u}+\lambda \mathbf{u}\in
\mathbf{E}^{s}(\Omega)$, если $\mathbf{u}\in
\mathbf{H}^{s+1}(\Omega)$.

Для функций $v$ из пространства $H^1(\Omega)$
определен \cite{mi} оператор {\it следа} $\gamma :H^1(\Omega)\rightarrow
H^{1/2}(\omega)$, равный следу $v$ на границе $\omega$ для  функций из $\mathcal{C}^1(\overline{\Omega})$: $\gamma v=v|_\omega$,
причем $\quad \|\gamma v\|_{L_2(\omega)}\leq c \|v\|_{H^1(\Omega)}$.

Аналогично, для  поля $\mathbf{u}(\mathbf{x})$ из
$\mathbf{E}^0(\Omega)$ определен \cite{rt} оператор {\it следа ее нормальной
	компоненты $\mathbf{n}\cdot \mathbf{u}$},  $\gamma_\mathbf{n} :\mathbf{E}^0(\Omega)\rightarrow
H^{-1/2}(\omega)$, равный сужению $\mathbf{n}\cdot \mathbf{u}$ на
$\omega$ для функций из $\mathcal{C}^1(\overline{\Omega})$:
$\gamma_{\mathbf{n}}
\mathbf{u}={\mathbf{n}}\cdot\mathbf{u}|_\omega$.

Для $\mathbf{u}\in \mathbf{E}^0(\Omega)$ и $v\in H^1(\Omega)$ верна
обобщенная формула Стокса:  \quad  $\langle \gamma_{\mathbf{n}} {\mathbf{u}},\gamma v \rangle=
(\mathbf{u}, \nabla v)+( \text{div}\mathbf{u}, v),$ \quad 
где $\langle\gamma_{\mathbf{n}} {\mathbf{u}},\gamma v \rangle$- линейный
функционал над пространством $H^{1/2}(\omega)$.  
  \quad 
   $\mathbf{{E}}_{\gamma}^0(\Omega)\equiv \{\mathbf{f}\in\mathbf{E}^0:\gamma_{\mathbf{n}}\mathbf{f}=0\}$.
\subsection {Метод Фурье  решения краевых задач в $\mathbf{L}_{2}(\Omega)$.
} Задано поле  $\mathbf{f}$ в  $\mathbf{L}_{2}(\Omega)$.

{\it Задача 1.   Найти поле 
	$\mathbf{u}$ в $	\mathbf{L}_{2}(\Omega)$ такое, что}
	\begin{equation*}	\klabel{rot__4_}
	 {\bf rot}\, \mathbf{u}+\lambda \mathbf{u}=
	\mathbf{f} \quad \text{в}\quad 	\mathbf{L}_{2}(\Omega),\eqno{(4.4)}\end{equation*}
то-есть  $(\mathbf{u}, ( {\bf rot}+\lambda\, I)\mathbf{ v})= (	\mathbf{f},	\mathbf{v})$ для любого поля $\mathbf{v}\in \,\mathbf{C}_0^{\infty}(\Omega).
 $

{\it Задача 2. Найти поле 
	$\mathbf{w}$ в $\mathbf{L}_{2}(\Omega)$ такое, что}
	\begin{equation*}	\klabel{nd_ 4_}
	 \nabla \mathbf {div}\, \mathbf{w}+\nu^2\mathbf{w}=
	\mathbf{f} \quad \text{в}\quad 	\mathbf{L}_{2}(\Omega).\eqno{(4.5)}\end{equation*}
    Переходим {\it  в объемлющее пространтво} $\mathbf{L}_{2}(\Omega)=  \mathcal{A} \oplus \mathcal{B}$. 
Используя разложение полей   $\mathbf{f}$ ,   $\mathbf{u}$ и  $\mathbf{w}$  
из $\mathbf{L}_2(\Omega)$ -- в  суммы:  $\mathbf{f}_\mathcal{A}+\mathbf{f}_\mathbf{V}$,   $\mathbf{u}_\mathcal{A}+\mathbf{u}_\mathbf{V}$  и $\mathbf{w}_\mathcal{A}+\mathbf{w}_\mathbf{V}$ и 
расширения  $S$ и $\mathcal{N}_d$ операторов ротор и  градиент дивергенции, эти уравнения запишем  в виде  уравнений-проекций :
	\begin{equation*}	\klabel{pro__4_}
	\lambda\mathbf{u}_\mathcal{A}=\mathbf{f}_\mathcal{A},  \quad (\mathcal{N}_d+\nu^2 I)\mathbf{w}_\mathcal{A}=\mathbf{f}_\mathcal{A}
	 \quad\text{в}\quad   \mathcal{A}
; \eqno{(4.6)}\end{equation*}
	\begin{equation*}	\klabel{prond__4_}
	 (S+\lambda I)\mathbf{u}_\mathbf{V}=\mathbf{f}_\mathbf{V},
\quad
\nu^2\mathbf{w}_\mathbf{V}=\mathbf{f}_\mathbf{V}, 
 \quad\text{в}\quad   \mathcal{B} \eqno{(4.7)}\end{equation*}
учитывая, что   $\text{rot} \,\mathbf{u}_\mathcal{A}=0$ в  $ \mathcal{A}$,   
$  \nabla\text{div}\,\mathbf{u}_\mathbf{V}=0$  в $ \mathcal{B}\equiv \mathbf{V}^0$.

{\it Замечание.} Если пространство    $\mathbf{C} \equiv \mathcal{B}_H(G)$ не пусто и     $\lambda\neq 0$,  то уравнение  $(\nabla\mathrm{div} +\mathrm{rot} +\lambda  I)\mathbf{u}=\mathbf{f}$ распадается на три проекции
\[   (\mathcal{N}_d+\lambda I)\mathbf{u}_\mathcal{A}=\mathbf{f}_\mathcal{A}, 
\quad (S+\lambda I)\mathbf{u}_\mathbf{V}=\mathbf{f}_\mathbf{V}, \quad  
\lambda\mathbf{u}_\mathbf{C}=\mathbf{f}_\mathbf{C}\eqno{(4.8)} \]
-  уравнения  второго, первого и нулевого порядков, соответственно.

 Согласно Теореме 3 и Леммам 1, 2, уравнения 
   $ (S+\lambda I)\mathbf{u}_\mathbf{V}=\mathbf{f}_\mathbf{V}$  и  $(\mathcal{N}_d+\nu^2 I)\mathbf{u}_\mathcal{A}=\mathbf{f}_\mathcal{A}$ разрешимы по Фредгольму. 

 При $\lambda\neq Sp (rot)$ проекции решения задачи 1 имеют вид:
	\begin{equation*} \klabel{kra 1}
{\mathbf{u}_{\mathcal{A}}}={\lambda}^{-1}\mathbf{f}_\mathcal{A}, \quad 
\mathbf{u}_\mathbf{V}=(S+\lambda)^{-1} \mathbf{f}_\mathbf{V}.
 \eqno{(4.9)}\end{equation*}

Действительно формулы  (4.9)  
 получаем из формул (4.6),   (4.7)   и обратимости оператора  $S+\lambda I$  при $\lambda\neq Sp (rot)$ в $ \mathbf{V}^0$ (Теорема 3 п.1.5).
   \quad Рассмотрим другие утверждения Теоремы 5: 
   
    если
 $\mathbf{f}=\mathbf{f}_{\mathcal{A}}$ и $\mathbf{f}_{\mathcal{A}}\in \mathcal{A}$ или  $\mathbf{f}_{\mathcal{A}}\in\mathcal{A}_{\gamma}$,
 то $\mathbf{u}=~{\lambda}^{-1}\mathbf{f}_{\mathcal{A}}\in
 \mathcal{A}$  или   $\mathbf{u}\in \mathcal{A}_{\gamma} $.
 
  Эти ряды являются обобщенными решениями уравнения (4.4).

 Если $\mathbf{f}\in \mathcal{B}\bot \mathcal{A}$ в
 $\mathbf{L}_{2}(\Omega)$, то
 $\mathbf{u}=(S+\lambda)^{-1}\mathbf{f}_\mathbf{V}\in \mathbf{W}^1\subset \mathbf{H}^{1}_{\gamma}(\Omega)$.
 
 	Если 	$\mathbf{f}\in\mathbf{E}_{\gamma}^{0}(\Omega)$,  то 	$\mathbf{u}=\mathbf{u}_{\mathcal{A}}+\mathbf{u}_\mathbf{V}\in  \mathbf{H}^{1}_{\gamma}(\Omega)$.
 	
 {\small 	Действительно, если $ \mathbf{f}\in
 	\mathbf{H}^{0}=\mathbf{L}_{2}(\Omega)$ и $ \text{div}\mathbf{f}\in
 	{H}^{0}$,  то $\mathbf{f}_{\mathcal{A}}\in 	\mathbf{H}^{0}$ и  $ \text{div}\,\mathbf{f}_{\mathcal{A}}=\text{div}\mathbf{f}\in
 	{H}^{0}$,  так как    $\text{div}\,\mathbf{f}_{\mathbf{V}}=0$.
 	Кроме того,    $\text{rot} \,\mathbf{f}_\mathcal{A}=0$ и $\gamma_{\mathbf{n}} {\mathbf{f}}_\mathcal{A}=0$. Из оценки (1.4) при $s=0$ поле
 	 $\mathbf{f}_{\mathcal{A}}\in 	\mathbf{H}^{1}_{\gamma}$, значит 
 	 $\mathbf{u}_{\mathcal{A}}={\lambda}^{-1}\mathbf{f}_\mathcal{A}\in \mathbf{H}^{1}_{\gamma}(\Omega)$. Так как $\mathbf{u}_\mathbf{V}=(S+\lambda)^{-1} \mathbf{f}_\mathbf{V}$ также принадлежит $	\mathbf{H}^{1}_{\gamma}$, то $\mathbf{u}\in  \mathbf{H}^{1}_{\gamma}(\Omega)$.}
 
 Если 	$\mathbf{f}\in C(2k, m)$,  то согласно	(4.9) $\mathbf{u}\in C(2k, m+1)$ .
 
 Если	же $\mathbf{f}\in\mathbf{C}_0^{\infty}(\Omega)$, то ряды (1.16),
 (1.17) сходятся в 	$\mathbf{H}^{s}(\Omega)$
 
 для 	любого $s\geq 1$ и $\mathbf{u}\in
 C^{\infty}(\overline{\Omega})$ - классическое решение задачи.

   Теорема 5 доказана.
   
   В случае шара эта теорема имеет наиболее натуральный  вид. 

{\small  Согласно \cite{saUMJ13} собственные значения $\lambda_{n,m}$  оператора $S$ в шаре радиуса $R$ равны $\pm \rho_{n,m}/R$,  где числа $\pm \rho_{n,m}$ - нули функций  $\psi_n(r)$,  (1.18),  $m,n\in {\mathbb {N}}$;  кратность  собственного значения  $\lambda_{n,m}$ равна $2n+1$. 
	Собственные поля 	 $\mathbf{q}^{\pm}_{\kappa}(\textbf{x})$	ротора  и   $\mathbf{q}_{\kappa}(\textbf{x})$ градиента дивергенци, $ \kappa=(n,m,k)$,  выражены явно через сферические функции и  функции $\psi_n(r)$ .}

	Из Теоремы 5 и Леммы 1  очевидно следуют Лемма 3 и Следствие.
	
 Решение краевой задачи 2	при $\lambda\neq Sp (\nabla \textbf{div})$  аналогично  \cite{ds18}.

Таким образом,  задачи 1 и 2 решены  полностью.

 \subsection{О приложениях}{\small  Собственные поля ротора
имеют приложения в
гидродинамике \cite{lad},
где они называются полями Бельтрами;
в астрофизике и в физике плазмы
они называются бессиловыми полями (force-free magnetic
fields - L.Woltjel \cite{wol},  free-decay fields -
J.B. Taylor \cite{tay}).
В.И.Арнольд \cite{ar} 1965 и В.В.Козлов \cite{koz} 1983
изучали топологию линий тока течений
идеальной жидкости при условии
$[\text{rot }\mathbf{v}, \mathbf{v}]\neq 0$.
Об этих и других  работах мы написали подробно в \cite{saUMJ13}.}  \quad

{\small  Отметим еще работы  L.Woltjel \cite{wol} и  \cite{wolc},  который ввел понятие спиральности (helicity) гладкого векторного поля в области $\Omega \subset \mathbb {R}^3$.

D. Cantarella, DeTec, G. Gluck and M.Tatel 2000 \cite{cdtgt}
  из лаборатории  "Физика плазмы" \, изучили линии тока
собственных полей ротора с минимальным собственным значением в шаре
и в шаровом слое.

Автор вывел формулы собственных полей ротора и градиента дивергенции в
шаре для любых собственных значений  (см.  \cite{sa2000}, \cite{sa01}, \cite{saUMJ13}). Формулы  собственных полей ротора, полученные   независимо от \cite{cdtgt} и опубликованные в \cite{sa2000} примерно в то же время,  дополняют формулы,  приведенные в \cite{cdtgt},  которые получены, следуя работам    \cite{wolc}\,  и  \cite{wol}.

Установлена связь собственных полей ротора и Стокса \cite{saVSTU}, \cite{saUMJ13}. Для нелинейной системы Навье-Стокса  с периодическими граничными условиями    найдены явные решения \cite{saTMF}. Совместно  с А.Г. Хайбуллиным  автор разработал
новый метод численного решения задачи Коши\cite{saha}, \cite{saUMJ}.

Профессор Исламов Г.Г.\cite{gais}, используя  формулу автора  из \cite{saVSTU} и  программу  Wolfram Mathematica  осуществил
визуализацию линий тока поля   $u_{1,1,0}^{+}(\mathbf {x})$ ротора с минимальным собственным значением в шаре ( см. http://www.wolfram.com/events/technology-conf.-ru/2016/
resources.html). Траектория движения  отдельной точки напоминает нить, которая 	наматывается на тороидальную катушку. 

В работе М.Е.Боговского \cite{bo} исследована задача Дирихле для оператора дивергенции, важная в гидродинамике (см. \cite{lad} $\S 2$ гл. 1).

Cтатья \cite{mabo}  В. Н.Масленниковой и М.Е. Боговского содержит обзор работ  по решению задачи С.Л.Соболева  \cite{sob54} и  \textit{аппроксимации потенциальных и соленоидальных векторных полей финитными бесконечно дифференцируемыми полями}.
В частности, они пишут, что в 1976 г. Heywood \cite{hey} построил соленоидальное векторное поле в $W_2^1(\Omega )$,  которое не аппроксимируется векторными полями из

 $J^{\infty}_{0}(\Omega)=\{\mathbf {u}\in C^{\infty}_{0}(\Omega) :\mathrm {div}\,\mathbf {u}=0 \}$.}

{\bf Благодарности: } Академику РАН, профессору В.П.Маслову,  профессору, доктору ф.-м. наук М.Д.Рамазанову
	 и  доценту, кандидату ф.-м. наук
Р.Н.Гарифуллину за  поддержку.

\newpage
Реферат:

Пространства Соболева  и операторы вихрь и градиент дивергенции

 Р.С.Сакс
 
  В ограниченной области $G$ с гладкой
 	границей изучаются 
 	краевые задачи для операторов ротор (вихрь) и градиент дивергенции с младшими   членами $\lambda \mathbf{u}$ в пространствах Соболева. 
 	
 	Особенность этих операторов состоит в том, что,  не являясь 
 	эллиптическими,  они  принадлежат при $\lambda\neq 0$  классу систем,  
 	приводимых  методом Б.Вайнберга и В.Грушина,  к эллиптическим матрицам.
 	 	И каждая из  задач удовлетворяет условиям эллиптичености
 	В.Солонникова.
 	Из элиптической теории и оценок вытекают  свойства решений спектральных задач ротора и
 	градиента дивергенции:
 	а) ненулевые собстенные значения имеют конечную кратность,
 	б)каждая  обобщенная собственная функция
 	бесконечно дифференцируема вплоть до границы области.
 	
Пространство  $\mathbf{L}_{2}(G)$ разлагается на два ортогональных подпространства потенциальных и соленоидальных полей,  которые  назовем классами и    обозначим через  $\mathcal{A}$ и  $\mathcal{B}$.  Они содержат подпространства $\mathcal{A}_{\gamma}\subset\mathcal{A}$  и  $\mathbf{V}^0\subset \mathcal{B}$. 
Операторы градиент дивергенции и ротор имеют самосопряженные
расширения $\mathcal{N}_d$ и $\mathcal{S}$ в 
подпространства $\mathcal{A}_{\gamma }$ и $\mathbf{V}^0$,  где они обратимы. Их обратные операторы  $\mathcal{N}_d^{-1}$ и $\mathcal{S}^{-1}$  -
вполне непрерывны, а  собственные векторы образуют ортогональные
базисы в каждом из классов  $\mathcal{A}$ и  $\mathcal{B}$. 
Их  элементы    $\mathbf{f}_\mathbf{A}$ и  $\mathbf{f}_\mathbf{V}$   разлагаются  
в  ряды Фурье,  а операторы  $\mathcal{N}_d: \mathcal{A}_{\gamma }\rightarrow\mathcal{A}_{\gamma } $ и $\mathcal{S}:\mathbf{V}^0  \rightarrow  \mathbf{V}^0$ 
получают спектральные представления. Их области определения $\mathcal{D}  (\mathcal{N}_d)\equiv \{\mathbf{f}\in \mathcal{A}_{\gamma }:\nabla  \mathrm{div} \mathbf{f}\in \mathcal{A}_{\gamma }\}$\quad  {и }  $\mathcal{D} (\mathcal{S})\equiv \{\mathbf{g}\in \mathbf{V}^0:\mathrm{rot}\mathbf{g}\in \mathbf{V}^0\}$
содержатся в пространствах Соболева   $\mathbf{H}^2$  и   $\mathbf{H}^1$   порядков 2 и 1. Мы выделяем  пространства
\[ \mathbf{A}^{2k}_{\gamma }\equiv \{\mathbf{f}\in \mathcal{A}_{\gamma },...,(\nabla  \mathrm{div})^k \mathbf{f}\in \mathcal{A}_{\gamma }\}\quad  \text{и } \mathbf{W}^m \equiv \{\mathbf{g}\in \mathbf{V}^0,..., (\mathrm{rot})^m \mathbf{g} \in \mathbf{V}^0\},\]  где  $k, m\geq 1$  и доказываем, что они являются аналогами пространств Соболева порядков  $2k$ и  $m$, соответственно,  в классах потенциальных и соленоидальных полей.
Прямые суммы этих пространств мы называем классами и обозначаем как $C(2k,m)$. 
Методом Фурье  при  $ \mathbf{f}\in C(2k,m)$ решаются краевые задачи 1 и 2 в  $\mathbf{L}_{2}(\Omega)$.   Доказано, что
 при  $\lambda\neq Sp(\mathrm{rot})$ оператор   $\text{rot}+\lambda I$ задачи 1 отображает  класс  $C(2k,m+1)$ на класс $C(2k,m)$ взаимно однозначно и непрерывно. \quad 
 Если $\lambda \overline{\in} Sp (\nabla \,\mathrm{div})$, 
 	то 	 оператор $\nabla \,\mathrm{div}+\lambda I$  задачи  2 отображает класс  $C(2(k+1), m)$
на  $C(2k, m)$	взаимно	однозначно и непрерывно. \quad 
В частности,  если область $\Omega=B$ есть шар,   $\psi_n(\lambda\,R)\neq 0$  для всех $n\in {\mathbb {N}}$ и функция  $\mathbf{f}\in  \mathbf{A}^{2k}_{\gamma }(B)\oplus  \mathbf{W}^m(B)$, то решение задачи 1 существует, единственно и принадлежит классу    $ \mathbf{A}^{2k}_{\gamma }(B)\oplus  \mathbf{W}^{m+1}(B)$. Соответственно, если   $\psi'_n(\nu\,R)\neq 0$  для всех $ 0\leq n\in {\mathbb {N}}$ и функция  $\mathbf{f}\in  \mathbf{A}^{2k}_{\gamma }(B)\oplus  \mathbf{W}^m(B)$, то решение задачи 2 существует, единственно и принадлежит классу    $ \mathbf{A}^{2(k+1)}_{\gamma }(B)\oplus  \mathbf{W}^{m}(B)$. 

Abstract: Saks Romen Semenovich

Sobolev spaces and operators vorticity and the gradient of the divergence

  In a bounded domain $G$ with smooth
 border studied 
 boundary value problems for operators of the rotor (vortex) and the gradient of the divergence with the  member  $\lambda \mathbf{u}$ in the Sobolev space.

The peculiarity of these operators is that they are not 
elliptic, they belong at $\lambda\neq 0$ to the class of systems,  
redused by the method of Veinberg and  Grushin  to elliptic matrices.
And each of these problems satisfies the  Solonnikov conditions of ellipticity.

From elliptic theory and estimates follow properties of solutions of spectral problems of the rotor and
 the gradient of the divergence:
a) non-zero eigenvalues have finite multiplicity,
b) each generalized eigenfunction be
infinitely differentiable up to the boundary of the domain.

It is known that the space $\mathbf{L}_{2}(G)$ decomposes into two orthogonal subspaces of potential and solenoidal fields, which we call classes and denote by $\mathcal{A}$ and $\mathcal{B}$.  They contain subspaces $\mathcal{A}_{\gamma}\subset\mathcal{A}$ and $\mathbf{V}^0\subset \mathcal{B}$.

It is proved that 
operators: gradient of the  divergence and rotor have self-adjoint
extensions $\mathcal{N}_d$ and $\mathcal{S}$ to  orthogonal
subspaces $\mathcal{A}_{\gamma }$ and $\mathbf{V}^0$, where they are reversible. 
Their inverse operators $\mathcal{N}_d^{-1}$ and $\mathcal{S}^{-1}$  
 are completely continuous, and their eigenvectors form an orthogonal
bases in each class $\mathcal{A}$ and $\mathcal{B}$.

The  elements  $\mathbf{f}_\mathbf{A}$ and $\mathbf{f}_\mathbf{V}$ decompose  
in Fourier series, and operators $\mathcal{N}_d: \mathcal{A}_{\gamma }\rightarrow\mathcal{A}_{\gamma } $ and $\mathcal{S}:\mathbf{V}^0 \rightarrow \mathbf{V}^0$ 
spectral representations are obtained. 
Their domains of definition is $\mathcal{D} (\mathcal{N}_d)\equiv \{\mathbf{f}\in \mathcal{A}_{\gamma }:\nabla \mathrm{div} \mathbf{f}\in \mathcal{A}_{\gamma}\} $ \quad {and} $\mathcal{D} (\mathcal{S})\equiv \{\mathbf{g}\in \mathbf{V}^0:\mathrm{rot}\mathbf{g}\in \mathbf{V}^0\}$
are contained in Sobolev spaces $\mathbf{H}^2$ and $\mathbf{H}^1$ of orders 2 and 1.

We introduce the spaces
\[ \mathbf{A}^{2k}_{\gamma }\equiv \{\mathbf{f}\in \mathcal{A}_{\gamma },..., (\nabla \mathrm{div})^k \mathbf{f}\in \mathcal{A}_{\gamma }\} \,\text{and } \mathbf{W}^m \equiv \{\mathbf{g}\in \mathbf{V}^0,..., \mathrm{rot}^m \mathbf{g} \in \mathbf{V}^0\}\] for $k, m\geq 1$  and prove they are analogues of the Sobolev space orders $2k$ and $m$, respectively, in the classes of potenti
al and solenoidal fields.
The direct sums of these spaces we call classes and denote as $C(2k,m)$.

The boundary value problems 1 and 2 in $\textbf{L}_2(\Omega)$ are solved by the Fourier method for $ \mathbf{f}\in C(2k,m)$.

If  $\lambda\neq Sp(\mathrm{rot})$, the $\textbf{rot}+\lambda I$ operator  maps class $C (2k, m+1)$ to class $C (2k,m)$ one-to-one and continuously.

  If $\lambda \overline{\in} Sp (\nabla \,\mathrm{div})$, 
 the  $\nabla \,\mathrm{div}+\lambda I$ operator  maps class $C(2(k+1), m)$ on   $C(2k, m)$ one-to-one and continuously.
 	
In particular,   if domain  $\Omega=B$ be a ball,   $\psi_n(\lambda\,R)\neq 0$ for all $n\in {\mathbb {N}}$ and field  $\mathbf{f}\in  \mathbf{A}^{2k}_{\gamma }(B)\oplus  \mathbf{W}^m(B)$, then solution of the problem  1 existe, unique and belongs to the class    $ \mathbf{A}^{2k}_{\gamma }(B)\oplus  \mathbf{W}^{m+1}(B)$.

 Respective, if    $\psi'_n(\nu\,R)\neq 0$  for all $ 0\leq n\in {\mathbb {N}}$ and the field \newline $\mathbf{f}\in  \mathbf{A}^{2k}_{\gamma }(B)\oplus  \mathbf{W}^m(B)$,  then solution of the problem  2 existe, unique and belongs to the class   $ \mathbf{A}^{2(k+1)}_{\gamma }(B)\oplus  \mathbf{W}^{m}(B)$.

Сакс Ромэн Семенович ведущий научный сотрудник
 Институт Математики с
ВЦ УФИЦ РАН 450077, г. Уфа, ул. Чернышевского, д.112 телефон: (347)
272-59-36
 (347) 273-34-12
факс: (347) 272-59-36 телефон дом.: (347) 273-84-69 моб.
+7 917 379 75 38

 e-mail: romen-saks@yandex.ru

 \end{document}